\documentclass[11pt,reqno]{amsart}

\usepackage[T1]{fontenc}
\usepackage[utf8]{inputenc}
\usepackage{lmodern}
\usepackage{microtype}
\usepackage{amsmath,amssymb,amsthm,mathtools,mathrsfs}
\usepackage{enumitem}
\usepackage{xcolor}
\usepackage{hyperref}
\usepackage[nameinlink,capitalise]{cleveref}
\usepackage[a4paper,margin=1.1in]{geometry}
\usepackage{comment}
\usepackage{tikz-cd}

\hypersetup{
  colorlinks=true,
  linkcolor=blue,
  citecolor=blue,
  urlcolor=blue,
  pdftitle={Strong Banach property (T) for SL3(R) and Banach spaces of nontrivial type}
}

\newtheorem{theorem}{Theorem}[section]
\newtheorem{proposition}[theorem]{Proposition}
\newtheorem{lemma}[theorem]{Lemma}
\newtheorem{corollary}[theorem]{Corollary}

\theoremstyle{definition}
\newtheorem{definition}[theorem]{Definition}
\newtheorem{remark}[theorem]{Remark}

\newcommand{\R}{\mathbb R}

\newcommand{\Z}{\mathbb Z}

\newcommand{\E}{\mathbb E}
\newcommand{\cF}{\mathcal F}
\newcommand{\cB}{\mathcal B}

\newcommand{\cL}{\mathcal L}

\newcommand{\cS}{\mathcal S}
\newcommand{\Id}{\mathrm{Id}}
\newcommand{\supp}{\operatorname{supp}}

\newcommand{\Rot}{\operatorname{Rot}}
\newcommand{\norm}[1]{\left\lVert #1\right\rVert}
\newcommand{\abs}[1]{\left\lvert #1\right\rvert}
\newcommand{\ip}[2]{\left\langle #1,#2\right\rangle}

\title[Strong Banach property $(T)$ for $\mathrm{SL}_3(\mathbb R)$]
{Strong Banach property $(T)$
for $\mathrm{SL}_3(\mathbb R)$}
\author{Jvbin Yao}
\address{Research Center for Operator Algebras,
	School of Mathematical Sciences,
	East China Normal University,
	Shanghai 200241, PR China}
\email{52285500009@stu.ecnu.edu.cn}

\begin{document}

\begin{abstract}
Strong Banach property~$(T)$ has been established for higher-rank
algebraic groups over non-Archimedean local fields with respect to
Banach spaces of nontrivial type. However, in the real setting,
existing results have required additional geometric hypotheses on the Banach
spaces, leaving the corresponding question open even for
$\mathrm{SL}_3(\mathbb R)$. We prove that $\mathrm{SL}_3(\mathbb R)$
has strong Banach property~$(T)$ with respect to every class of
complex Banach spaces of uniformly nontrivial Rademacher type.
As a consequence, the same holds for every lattice in
$\mathrm{SL}_3(\mathbb R)$, in particular for $\mathrm{SL}_3(\mathbb Z)$.
\end{abstract}

\maketitle

\section{Introduction}

For a topological group $G$, Kazhdan's property~$(T)$
is a rigidity property for unitary representations of $G$
(see \cite{BHV}). Its connections with fixed-point properties for affine
isometric actions on Hilbert spaces and with expander graphs
have made it a fundamental tool in rigidity theory and
geometric group theory
(see \cite[Theorem 2.12.4 and Corollary 6.1.9]{BHV}).
Rigidity problems for group actions on manifolds provide
a natural motivation for extending property~$(T)$ beyond
the Hilbert space setting. In particular, Fisher and
Margulis showed that $L^p$-analogues of property~$(T)$
for sufficiently large $p$ can be used to weaken
smoothness assumptions in certain rigidity results
\cite{FM1,FM2}; see also \cite[Section~1.b]{BFGM}.
This was one of the motivations for the systematic study
of Banach space analogues of property~$(T)$ by Bader,
Furman, Gelander, and Monod \cite{BFGM}, who also 
investigated their connections with fixed-point properties for affine 
isometric actions.

At roughly the same time as the work of \cite{BFGM}, in the course
of investigating the Baum--Connes conjecture with
arbitrary coefficients using Banach $KK$-theory
\cite{LafforgueBanachKK,LafforgueBanachKKInvent}, Lafforgue introduced a
strengthening of Kazhdan's property~$(T)$ in
\cite{LafforgueStrongT}, now known as strong
property~$(T)$. 
His original work already contained a Banach space version of this strengthening for certain classes of Banach spaces 
defined by Fourier transform estimates (see \cite[Definition~0.2]{LafforgueStrongT}). 
Lafforgue \cite[Introduction]{LafforgueStrongT} noted that the question of extending
strong property~$(T)$ to Banach spaces  had been raised by Uri Bader in connection with
the work of Bader, Furman, Gelander, and Monod \cite[Remark~1.7(2)]{BFGM}.
In a subsequent paper, he strengthened his results in the non-Archimedean setting by establishing 
strong Banach property~$(T)$ with respect to classes of Banach spaces of uniformly nontrivial Rademacher type
\cite[Definitions~0.1--0.2 and Theorem~0.3]{LafforgueBanachT}(see \cref{def:rademacher-type} for the definition of
Rademacher type).
These results also yield fixed-point properties for affine isometric actions on Banach spaces of nontrivial type and
families of expanders that do not admit uniform embeddings into any Banach space of nontrivial type
\cite[Corollaries~0.5 and~0.7]{LafforgueBanachT}.

The fixed-point and nonembeddability consequences of Lafforgue's results make it natural to ask how 
broadly strong Banach property~$(T)$ holds, both in terms of the groups and the Banach spaces involved. 
In the non-Archimedean setting, Lafforgue's theorem applies to $\mathrm{SL}_3(F)$ and, more generally, 
to connected almost $F$-simple algebraic groups whose Lie algebras contain $\mathfrak{sl}_3(F)$, where 
$F$ is a non-Archimedean local field \cite{LafforgueBanachT}. Liao subsequently extended this result to 
$\mathrm{Sp}_4(F)$ and thereby to all connected almost $F$-simple algebraic groups of $F$-split rank at 
least two \cite{LiaoStrongBanachT}. Thus, in the non-Archimedean higher-rank setting, strong Banach property~$(T)$ 
is known for every class of Banach spaces of uniformly nontrivial Rademacher type.

The results over non-Archimedean local fields naturally raise the question of whether analogous results hold 
over Archimedean local fields, particularly over $\mathbb R$. For $\mathrm{SL}_3(\mathbb R)$, Lafforgue first 
established strong property~$(T)$ with respect to Hilbert spaces \cite{LafforgueStrongT}. De la Salle \cite{deLaSalleSL3} 
subsequently extended this result to the class $\mathcal E_4$ of Banach spaces. The class $\mathcal E_4$ is defined using 
type and cotype conditions, complex interpolation, and certain closure operations (see \cite[Definition~3.2]{deLaSalleSL3}). 
It contains all $L^p$-spaces for $1<p<\infty$, as well as some nonreflexive Banach spaces (see \cite[Sections~3.1--3.2]{deLaSalleSL3}).

The extension from $\mathrm{SL}_3(\mathbb R)$ to general higher-rank real Lie groups was subsequently achieved by de Laat and 
de la Salle \cite{deLaatDeLaSalleStrongT}. They proved strong Banach property~$(T)$ for every connected simple Lie group of 
real rank at least two with respect to the class $\mathcal E_{10}\subseteq\mathcal E_4$ (see \cite[Theorem~A]{deLaatDeLaSalleStrongT}). 
In another direction, de Laat, Mimura, and de la Salle \cite[Theorem~A]{deLaatMimuraDeLaSalle} showed that more general geometric conditions 
on Banach spaces can be allowed when the real rank is sufficiently large.

The extension of strong Banach property~$(T)$ to lattices presents further difficulties. Earlier results also covered
cocompact lattices by induction of representations, but this approach does not readily extend to non-cocompact lattices,
such as $\mathrm{SL}_3(\Z)$. De la Salle \cite{deLaSalleLattices} overcame this difficulty by introducing two-step representations and a stronger
property~$(\ast_{\mathcal E})$. He established strong Banach property~$(T)$ for arbitrary lattices in higher-rank groups
over non-Archimedean local fields with respect to classes of Banach spaces of uniformly nontrivial type, and obtained
analogous results in the real setting under additional assumptions on the Banach spaces \cite[Theorems~1.4--1.5]{deLaSalleLattices}.

Despite these developments, nontrivial Rademacher type alone was not known to suffice for strong Banach property~$(T)$
even for $\mathrm{SL}_3(\R)$. In view of the non-Archimedean results, 
Lafforgue asked whether strong Banach property~$(T)$ holds for higher-rank algebraic
groups over Archimedean local fields, in particular for $\mathrm{SL}_3(\R)$ \cite[Introduction, question~(b)]{LafforgueBanachT}
(see also \cite[Introduction, question~(b)]{LiaoFFT}).

In this paper, we answer Lafforgue's question affirmatively for $\mathrm{SL}_3(\R)$ and obtain the corresponding result
for all its lattices. To state it precisely, we recall the following formulation of strong Banach property~$(T)$ 
from \cite[Definition~1.1]{deLaSalleSL3}.

\begin{definition}\label{def:strong-banach-T}
Let $G$ be a locally compact group and $\mathcal E$ a class of 
complex Banach spaces. We say that $G$ has strong Banach property~$(T)$ 
with respect to $\mathcal E$ if, for every length function $\ell$ on $G$, 
there exists a sequence $(m_n)$ of compactly supported symmetric Borel measures 
on $G$ with the following property.

For every $E\in\mathcal E$, there exists $t>0$ such that
every strongly continuous representation $\pi$ of $G$ on $E$
satisfying $\|\pi(g)\|_{\mathcal B(E)}\le Le^{t\ell(g)}$
for all $g\in G$ and some $L>0$
admits a bounded projection $P_\pi$ onto $E^{\pi(G)}$ such that
\[
\|\pi(m_n)-P_\pi\|_{\mathcal B(E)}\longrightarrow0.
\]
Here
$E^{\pi(G)}=\{v\in E:\pi(g)v=v\text{ for all }g\in G\}$
and $\pi(m_n)v=\int_G\pi(g)v\,dm_n(g)$.
Symmetry means $m_n(A^{-1})=m_n(A)$
for every Borel set $A\subset G$.

\end{definition}

A class $\mathcal E$ of complex Banach spaces is said to 
have type $>1$ if there exist $p>1$ and $T<\infty$ such that 
every $E\in\mathcal E$ has Rademacher type $p$ with constant 
at most $T$ (see \cref{def:rademacher-type}).

With this terminology, our main theorem is as follows.

\begin{theorem}\label{thm:main}
$\mathrm{SL}_3(\mathbb R)$ has strong Banach property $(T)$
with respect to every class $\mathcal E$ of complex
Banach spaces of type $>1$.
\end{theorem}

Using the same analytic estimates together with de la Salle's
criterion for lattices \cite[Theorem~1.5]{deLaSalleLattices},
we also obtain the following result.

\begin{corollary}\label{cor:lattice-strongT}
Let $\mathcal E$ be a class of complex Banach spaces of type $>1$.
Every lattice $\Gamma$ in $\mathrm{SL}_3(\R)$ has strong Banach
property~$(T)$ with respect to $\mathcal E$, with the growth
exponent in \cref{def:strong-banach-T} uniform over $\mathcal E$
for each length function on $\Gamma$.
In particular, this holds for $\mathrm{SL}_3(\Z)$.
\end{corollary}

For $x\in S^2$ and $|\delta|<1$, let $(T_\delta f)(x)$
denote the average of $f$ over the circle
$\{y\in S^2:\langle x,y\rangle=\delta\}$. At the endpoints,
set $T_1f(x)=f(x)$ and $T_{-1}f(x)=f(-x)$
(see \cref{sec:sphere} for details).
By the criterion due to Lafforgue
\cite[Introduction, question~(b)]{LafforgueBanachT}
and proved in detail by de la Salle
\cite[Theorem~4.1]{deLaSalleSL3},
together with de la Salle's criterion for lattices
\cite[Theorem~1.5]{deLaSalleLattices},
both \cref{thm:main} and \cref{cor:lattice-strongT}
reduce to establishing the following estimate:
there exist constants $C,\alpha>0$ such that
\[
\sup_{E\in\mathcal E}
\bigl\|(T_\delta-T_0)\otimes\Id_E
\bigr\|_{\mathcal B(L^2(S^2;E))}
\le C|\delta|^\alpha,
\qquad \delta\in[-1,1].
\]

De la Salle established a corresponding estimate for each 
$E\in\mathcal E_4$. His approach combines scalar Schatten estimates of order 
$r>4$ with a vector-valued bound requiring $1/p-1/q<1/r$ for spaces of type $p$ 
and cotype $q$. Thus the direct argument requires $1/p-1/q<1/4$, a condition not 
guaranteed by nontrivial type alone. Although complex interpolation and closure 
operations enlarge the class of spaces covered, the resulting arguments do not yield 
the estimate for arbitrary Banach spaces of nontrivial type \cite[Lemma~2.13 and Section~3]{deLaSalleSL3}. 
This motivates the search for an approach that establishes the required decay under a common 
nontrivial type bound alone.

Our approach is to exploit the decay of finite Fourier transforms directly, rather than to deduce
vector-valued estimates from scalar Schatten estimates. The main analytic step is to convert finite 
Fourier transform decay into decay for integral operators satisfying a rotational curvature condition.

To see how such integral operators arise, we approximate the
spherical averaging operators $T_\delta$ by integral operators
with smooth kernels. Using an approximate identity and a
telescoping decomposition, we write
$T_\delta=\mathsf S_\delta+\sum_{j=0}^{\infty}\mathsf D_{\delta,j}$.
Both $\mathsf S_\delta$ and $\mathsf D_{\delta,j}$ have kernels
of the form $Rh(R\rho_\delta(x,y))$, where
$\rho_\delta(x,y)=\langle x,y\rangle-\delta$.
The key distinction is that the function $h$ associated with
$\mathsf D_{\delta,j}$ has Fourier transform supported away
from zero. The required H\"older estimate then reduces to
establishing uniform power decay for these difference operators
(see \cref{sec:sphere}). This leads us to study the decay of integral operators
of the more general form
\[
    (\mathcal K_R f)(q)
    =
    \int_P Rk(R\rho(q,p))a(q,p)f(p)\,d\mu_P(p),
\]
where $P$ and $Q$ are two-dimensional smooth manifolds
equipped with smooth positive densities,
$a\in C_c^\infty(Q\times P)$, and $k\in\mathcal S(\R)$ with
$\operatorname{supp}\widehat{k}\Subset\R\setminus\{0\}$.
The functions $\rho_\delta(x,y)=\langle x,y\rangle-\delta$
satisfy a uniform nondegenerate rotational curvature condition
for $|\delta|\le 1/2$ (It is unclear to us how far this curvature condition can be weakened).
This suggests establishing a decay estimate based on
rotational curvature.

In fact, we prove in \cref{thm:curved} that the condition
$\operatorname{Rot}(\rho)\ne0$ on
$\operatorname{supp}a\cap\{\rho=0\}$ is sufficient to obtain
\[
    \sup_{E\in\mathcal E}
    \|\mathcal K_R\otimes\Id_E\|
    \le CR^{-\gamma},
    \qquad R\ge1,
\]
for some $\gamma>0$, whenever $\mathcal E$ is a class
of type $>1$ and $\rho$ is sufficiently smooth. Although we do not pursue 
other applications here, the resulting theorem provides a general Banach-valued
rotational curvature estimate under a nontrivial Rademacher type assumption alone.

To establish this estimate, we first reduce the general
integral operators to model operators by an anisotropic
localization. By the rotational curvature condition, we can write
in suitable product coordinates
\[
    H^{-1}\rho(q,p)
    =t-t'+cxy+L\Theta_L(y,t,x,t'),
    \qquad H=L^2,\quad L=R^{-1/3},
\]
where $|c|$ is bounded away from zero and $\Theta_L$ is
uniformly bounded in a suitable $C^r$ norm (see \cref{lem:plate-cover}\textup{(ii)}).
Since $RL^3=1$, the remainder can be absorbed into the
amplitude after Fourier inversion.
A Fourier series decomposition of the resulting amplitude
then reduces the problem to estimating model operators
of the form
\[
    (\mathsf S_{M,c,\psi}F)(y,t)
    =
    \iint M\check{\psi}\bigl(M(t-t'+cxy)\bigr)
    F(x,t')\,dt'\,dx,
\]
where $M=RH=R^{1/3}$ and $\psi$ is supported away from zero
(see \cref{cor:flat-amplitude}).

To estimate $\mathsf S_{M,c,\psi}$, we restrict the Fourier
transform in the $t$-variable to a short interval centered
at $\xi_0$, with $|\xi_0|\asymp M$. Removing the central
frequency leads to operators of the form
\[
    (\mathsf A_{\xi_0,c}g)(y,t)=\int_0^1 e^{ic\xi_0xy}g(x,t+cxy)\,dx.
\]
This suggests using the finite Fourier factorization
method developed in the non-Archimedean setting by
Lafforgue \cite[Lemma~3.2]{LafforgueBanachT} and
Liao \cite[Proposition~3.1 and Corollary~3.2]{LiaoFFT}.
The main difficulty here is to incorporate the continuous
translations in $t$ into the factorization without losing
the decay of the finite Fourier transforms.

After suitable changes of variables, the problem reduces
to estimating operators of the form
\[
    (\mathcal T_{N,\beta}g)(y,t)
    =
    \int_0^1 e^{2\pi iA^Nxy}g(x,t+\beta xy)\,dx,
\]
where $A\ge2$ is an integer, $\beta\in[0,1]$, and
$A^N\asymp M$. Expanding $x$ and $y$ in base $A$, we decompose
$xy=h_0+h_1$, where $A^Nh_0$ is integer-valued.
Although $h_0$ makes no contribution to the oscillatory
factor, it still contributes to the translation in $t$.
The key observation is that the dependence of $h_0$ on
the digits allows the corresponding translations to be
inserted as isometries between successive finite Fourier
transforms. This gives a factorization
\[
    S_N=JF_ND_N\cdots F_1D_1,
\]
where $\|J\|\le1$, $\|D_j\|=1$, and
$\|F_j\|\le\rho_A(E)$, with
$\rho_A(E)=\|\mathcal F_A\otimes\Id_E\|$.
Consequently, $\|S_N\|\le\rho_A(E)^N$. The remaining phase factors
and translations are controlled by scalar and operator-valued Schur multiplier estimates.
For $g$ with Fourier support in $[-\Delta,\Delta]$ in the
$t$-variable, we obtain
\[
    \|\mathcal T_{N,\beta}g\|
    \le C_0
    \bigl(C_1e^{C_2\vartheta A}\rho_A(E)\bigr)^N
    \|g\|,
\]
provided that $\beta\Delta\le 2\pi\vartheta A^N$
(see \cref{prop:digit}). By nontrivial type, $\sup_{E\in\mathcal E}\rho_A(E)\le CA^{-\eta}$
for some $\eta>0$. Combined with the preceding estimate, this yields
exponential decay in $N$, and hence power decay in $M$ for $\mathsf S_{M,c,\psi}$.
Combining the local estimates yields the desired decay
for $\mathcal K_R$. Applying this estimate to the spherical
averaging operators gives the required H\"older estimate.

The rest of the paper is organized as follows.
\Cref{sec:prelim} recalls a finite Fourier transform estimate and
states the criterion for strong Banach property~$(T)$.
\Cref{sec:flat} develops the finite Fourier factorization
in the continuous setting and establishes power decay for
the model operators.
\Cref{sec:curved} extends this decay estimate to integral
operators satisfying a rotational curvature condition.
\Cref{sec:sphere} uses this estimate to prove the H\"older
estimate for spherical averaging operators, and
\cref{sec:strongT} completes the proofs of the main theorem
and its lattice consequence.

\section{Preliminaries}
\label{sec:prelim}
In this section, we present a decay estimate for finite
Fourier transforms and a criterion for strong Banach
property $(T)$ in terms of operator norm estimates
for spherical averages.
The former is used in \cref{sec:flat}, and the latter
in the proof of \cref{thm:main}.

We work throughout with complex Banach spaces.
We begin by recalling the definition of Rademacher type.

\begin{definition}\label{def:rademacher-type}
A Banach space $E$ has Rademacher type $p\in[1,2]$ if there is $T<\infty$
such that for every finite family $x_1,\dots,x_n\in E$,
\[
 \left(\E_\varepsilon
 \norm{\sum_{j=1}^n\varepsilon_jx_j}^{2}\right)^{1/2}
 \le T\left(\sum_{j=1}^n\norm{x_j}^{p}\right)^{1/p}.
\]
Here $\varepsilon_1,\ldots,\varepsilon_n$ are independent
random variables taking the values $\pm1$
with equal probability.
The least admissible constant $T$ is denoted by $T_p(E)$,
with $T_p(E)=\infty$ if no finite constant exists.
The space $E$ has nontrivial type if this holds
for some $p>1$.
\end{definition}

For classes of type $>1$, the norms of normalized
finite Fourier transforms decay uniformly as the
group order increases.
To state this estimate precisely, we introduce
the following notation.

For an integer $A\ge2$, equip
$G_A=\mathbb Z/A\mathbb Z$ with normalized counting
measure and write
$\ell_A^2(E)=L^2(G_A;E)$ and
$\ell_A^2=\ell_A^2(\mathbb C)$.
Identifying the dual group of $G_A$ with $G_A$
via the characters $a\mapsto e^{2\pi iab/A}$,
define the normalized finite Fourier transform
$\mathcal F_A:\ell_A^2\to\ell_A^2$ by
\[
(\mathcal F_Af)(b)
=\frac1A\sum_{a\in G_A}e^{2\pi iab/A}f(a),
\qquad b\in G_A.
\]
Using the representatives $0,\ldots,A-1$ of
$\mathbb Z/A\mathbb Z$, we shall henceforth write
\[
(\cF_Af)(b)
=\frac1A\sum_{a=0}^{A-1}e^{2\pi iab/A}f(a),
\qquad b=0,\ldots,A-1.
\]
For a Banach space $E$, set
$\rho_A(E)
=
\|\mathcal F_A\otimes\Id_E\|_{\ell_A^2(E)\to\ell_A^2(E)}.$
If $E$ is a nonzero Hilbert space, then
$\rho_A(E)=A^{-1/2}$.


The following consequence of Bourgain's theorem
\cite[Theorem~1]{BourgainHY} is stated in
\cite[Corollary~2.2]{LiaoFFT}.

\begin{theorem}[Liao {\cite[Corollary~2.2]{LiaoFFT}}]
\label{thm:liao}
Let $\mathcal E$ be a class of type $>1$.
There exist constants $C>0$ and $\eta>0$ such that
\[
\rho_A(E)\le CA^{-\eta}
\qquad(A\ge2,\ E\in\mathcal E).
\]
\end{theorem}

We will also need finite Fourier estimates for functions
taking values in closed subspaces of $E$-valued $L^2$ spaces.

\begin{lemma}\label{lem:bochner}
Let $E$ be a Banach space, $\Omega$ a measure space,
and $H$ a closed subspace of $L^2(\Omega;E)$.
Then $\rho_A(H)\le\rho_A(E)$ for every integer $A\ge2$.
If $E$ has Rademacher type $p\in[1,2]$, then
$T_p(H)\le T_p(E)$.
\end{lemma}

\begin{proof}
	We first prove both estimates for $L^2(\Omega;E)$.
	Let $J:\ell_A^2(L^2(\Omega;E))
	\longrightarrow
	L^2(\Omega;\ell_A^2(E))$
	be the canonical identification $J((f_a)_{a=0}^{A-1})(\omega)=(f_a(\omega))_{a=0}^{A-1}$.
	By Fubini, \(J\) is an isometric isomorphism.
	For \(F=(f_a)_{a=0}^{A-1}\in \ell_A^2(L^2(\Omega;E))\) and for almost every
	\(\omega\in\Omega\), we have 
	\[
	J\bigl((\mathcal F_A\otimes\Id_{L^2(\Omega;E)})F\bigr)(\omega)
	=
	(\mathcal F_A\otimes\Id_E)\bigl(JF(\omega)\bigr),
	\]
	since the \(b\)-th coordinate of both sides is
	$
	\frac1A\sum_{a=0}^{A-1}e^{2\pi iab/A}f_a(\omega).
	$
	Hence, by the definition of \(\rho_A(E)\),
	\[
	\left\|
	J\bigl((\mathcal F_A\otimes\Id_{L^2(\Omega;E)})F\bigr)(\omega)
	\right\|_{\ell_A^2(E)}
	\le
	\rho_A(E)\,\|JF(\omega)\|_{\ell_A^2(E)}
	\]
	for almost every \(\omega\).  Integrating the squares over \(\Omega\) and
	using that \(J\) is an isometry gives
	\[
	\|(\mathcal F_A\otimes\Id_{L^2(\Omega;E)})F\|_{\ell_A^2(L^2(\Omega;E))}
	\le
	\rho_A(E)\,
	\|F\|_{\ell_A^2(L^2(\Omega;E))}.
	\]
	Thus $\rho_A(L^2(\Omega;E))\le \rho_A(E)$.
	
	Now suppose that \(E\) has Rademacher type \(p\) with constant \(T\), and
	let \(f_1,\dots,f_n\in L^2(\Omega;E)\).  By Fubini and the type-\(p\)
	inequality in \(E\),
	\[
	\left(
	\E_\varepsilon
	\left\|\sum_{j=1}^n\varepsilon_jf_j\right\|_{L^2(\Omega;E)}^2
	\right)^{1/2}
	\le
	T
	\left(
	\int_\Omega
	\left(\sum_{j=1}^n\|f_j(\omega)\|_E^p\right)^{2/p}
	\,d\omega
	\right)^{1/2}.
	\]
	Since \(p\le2\), we have \(2/p\ge1\), and the triangle inequality in
	\(L^{2/p}(\Omega)\) gives
	\[
	\left(
	\int_\Omega
	\left(\sum_{j=1}^n\|f_j(\omega)\|_E^p\right)^{2/p}
	\,d\omega
	\right)^{1/2}
	\le
	\left(
	\sum_{j=1}^n
	\|f_j\|_{L^2(\Omega;E)}^p
	\right)^{1/p}.
	\]
	Therefore \(L^2(\Omega;E)\) has Rademacher type \(p\), with type constant
	at most \(T\).
	Thus $T_p(L^2(\Omega;E))\le T_p(E)$.
	Since $H$ carries the inherited norm, restricting the finite
	Fourier transform to $\ell_A^2(H)$ and the type inequality
	to finite families in $H$ gives both conclusions.
\end{proof}

To formulate the strong Banach property $(T)$ criterion
in terms of spherical averages, we first compare
the averaging operators on $\mathrm{SO}(3)$ and $S^2$.
Let \(K=\mathrm{SO}(3)\), \(e_1=(1,0,0)\), and let
\[
U=\operatorname{Stab}_K(e_1)
=
\left\{
\begin{pmatrix}
	1&0\\
	0&A
\end{pmatrix}
:A\in\mathrm{SO}(2)
\right\}
\simeq\mathrm{SO}(2).
\]
Let \(\lambda\) be the left regular representation of \(K\) on \(L^2(K)\). As recalled in \cite{deLaSalleSL3}, the double coset space \(U\backslash K/U\) is parametrized by \(UxU\longmapsto x_{11}\in[-1,1]\). Thus, for \(\delta\in[-1,1]\), if \(x_\delta\in K\) satisfies \((x_\delta)_{11}=\delta\), we set \(T_\delta^K=\lambda(Ux_\delta U)\).

The operator \(T_\delta^K\) preserves \(L^2(U\backslash K)\) and
vanishes on its orthogonal complement. Under the identification
$U\backslash K\simeq S^2$ with $Uk\longmapsto k^{-1}e_1$,
its restriction to \(L^2(U\backslash K)\) corresponds to the spherical
averaging operator \(T_\delta\) given, for continuous \(f\), by
\[
(T_\delta^{S^2}f)(z)
=
\int_{\{y\in S^2:\langle z,y\rangle=\delta\}}
f(y)\,d\sigma_{z,\delta}(y),
\]
where \(\sigma_{z,\delta}\) denotes the normalized rotation-invariant
probability measure on this circle; see
\cite[\S2.10, (3) and the discussion immediately following it]
{deLaSalleSL3}. For later use, we establish the corresponding norm
identity on Banach-valued $L^2$ spaces.

\begin{lemma}\label{lem:K-sphere-comparison}
	For every Banach space \(E\) and every \(\delta\in[-1,1]\),
	\[
	\bigl\|
	(T_\delta^K-T_0^K)\otimes\Id_E
	\bigr\|_{L^2(K;E)\to L^2(K;E)}
	=
	\bigl\|
	(T_\delta^{S^2}-T_0^{S^2})\otimes\Id_E
	\bigr\|_{L^2(S^2;E)\to L^2(S^2;E)}.
	\]
\end{lemma}

\begin{proof}
	Let \(P_U:L^2(K;E)\to L^2(K;E)\) be the left \(U\)-averaging operator
	$(P_Uf)(k)=\int_U f(u^{-1}k)\,du$.
	Then \(P_U\) is a contractive projection, and
	\[
	\operatorname{Ran}P_U
	=
	\bigl\{
	f\in L^2(K;E):
	f(uk)=f(k)
	\text{ for every }u\in U
	\text{ and a.e. }k\in K
	\bigr\}.
	\]
	Hence $\operatorname{Ran}P_U$ identifies isometrically
    with $L^2(U\backslash K;E)$, equipped with the quotient
    probability measure.
	
	Since \(T_\delta^K=P_U\lambda(x_\delta)P_U\)  in the scalar case,
	its \(E\)-valued extension satisfies $T_\delta^K\otimes\Id_E
	=
	P_U\bigl(\lambda(x_\delta)\otimes\Id_E\bigr)P_U$.
	Set \(A_\delta=(T_\delta^K-T_0^K)\otimes\Id_E\). Then
	\(A_\delta=A_\delta P_U\) and \(P_UA_\delta=A_\delta\). Thus \(A_\delta|_{\operatorname{Ran}P_U}\) defines an operator on
	\(\operatorname{Ran}P_U\).
	
	We claim that $\|A_\delta\|=\bigl\|A_\delta|_{\operatorname{Ran}P_U}\bigr\|$.
	Indeed, the inequality
	\(\bigl\|A_\delta|_{\operatorname{Ran}P_U}\bigr\|\le \|A_\delta\|\)
	is immediate. Conversely, since \(A_\delta=A_\delta P_U\) and \(P_U\)
	is contractive, for every \(f\in L^2(K;E)\),
	\[
	\|A_\delta f\|
	=
	\|A_\delta P_Uf\|
	\le
	\bigl\|A_\delta|_{\operatorname{Ran}P_U}\bigr\|
	\,\|P_Uf\|
	\le
	\bigl\|A_\delta|_{\operatorname{Ran}P_U}\bigr\|
	\,\|f\|.
	\]
	Thus the reverse inequality also holds.
	
	Let $d\sigma$ denote surface area measure on $S^2$.
    Under $Uk\mapsto k^{-1}e_1$, the quotient probability
    measure corresponds to $d\sigma/(4\pi)$. 
    Define $J:\operatorname{Ran}P_U\to L^2(S^2,d\sigma;E)$ by
    $(Jf)(k^{-1}e_1)=\frac{1}{\sqrt{4\pi}}f(k)$ for almost every $k\in K$.
    Then $J$ is an isometric isomorphism.
	By the spherical realization recalled above, after tensoring
	with \(\Id_E\), the following diagram commutes:
   \[
   \begin{tikzcd}[column sep=7em, row sep=3.5em]
   	\operatorname{Ran}P_U
   	\arrow[r, "{A_\delta|_{\operatorname{Ran}P_U}}"]
   	\arrow[d, "J"']
   	&
   	\operatorname{Ran}P_U
   	\arrow[d, "J"]
   	\\
   	L^2(S^2;E)
   	\arrow[r, "{(T_\delta^{S^2}-T_0^{S^2})\otimes\Id_E}"']
   	&
   	L^2(S^2;E)
   \end{tikzcd}
   \]
	Hence $J\bigl(A_\delta|_{\operatorname{Ran}P_U}\bigr)J^{-1}=
	(T_\delta^{S^2}-T_0^{S^2})\otimes\Id_E$. 
	Since \(J\) is an isometric isomorphism,
	$\bigl\|A_\delta|_{\operatorname{Ran}P_U}\bigr\|
	=
	\bigl\|
	(T_\delta^{S^2}-T_0^{S^2})\otimes\Id_E
	\bigr\|_{L^2(S^2;E)\to L^2(S^2;E)}$.
	Combining this with the preceding equality of norms gives
	\[
	\bigl\|
	(T_\delta^K-T_0^K)\otimes\Id_E
	\bigr\|_{L^2(K;E)\to L^2(K;E)}
	=
	\bigl\|
	(T_\delta^{S^2}-T_0^{S^2})\otimes\Id_E
	\bigr\|_{L^2(S^2;E)\to L^2(S^2;E)}.
	\]
\end{proof}

In the following criterion, write $T_\delta=T_\delta^{S^2}$.
By \cref{lem:K-sphere-comparison}, the spherical estimate
below is equivalent to Assumption~(6) of \cite{deLaSalleSL3},
with $s=\alpha$.
The criterion therefore follows from
\cite[Theorem~4.1 and the paragraph immediately following it]
{deLaSalleSL3}.

\begin{theorem}[Lafforgue--de la Salle
{\cite[Theorem~4.1]{deLaSalleSL3}}]
\label{thm:criterion}
Let $\mathcal E$ be a class of complex Banach spaces.
Suppose there exist $C>0$ and $\alpha\in(0,1/2]$ such that
\[
\norm{(T_\delta-T_0)\otimes\Id_E}_{L^2(S^2;E)\to L^2(S^2;E)}
\le C|\delta|^\alpha
\qquad(E\in\mathcal E,\ \delta\in[-1,1]).
\]
Then $\mathrm{SL}_3(\R)$ has strong Banach property $(T)$
with respect to $\mathcal E$.
\end{theorem}

\begin{remark}\label{rem:criterion-uniform}
In \cref{thm:criterion}, for each length function $\ell$,
the exponent $t$ in \cref{def:strong-banach-T}
can be chosen to depend only on $C,\alpha,\ell$,
independently of $E\in\mathcal E$.
This follows from \cite[Theorem~4.1]{deLaSalleSL3}
and the bound $\ell(g)\le a\ell_0(g)+b$ with $g\in\mathrm{SL}_3(\R)$,
where $\ell_0(g)=\max\{\log\|g\|,\log\|g^{-1}\|\}$,
with the Euclidean operator norm, and $a>0$, $b\ge0$
depend only on $\ell$; see \cite[Section~4.1]{deLaSalleSL3}.
\end{remark}

\section{Integral operator estimates via finite Fourier transforms}
\label{sec:flat}

\subsection{Schur multipliers induced by projective tensor products}

Let $(U,\mu_U)$ and $(V,\mu_V)$ be probability spaces,
and let $H$ be a Banach space. Set
\[
\mathfrak A(V,U)
=L^\infty(V)\widehat\otimes_\pi L^\infty(U),
\qquad
\mathfrak B(V,U;H)
=L^\infty(V)\widehat\otimes_\pi
\mathcal L(H)\widehat\otimes_\pi L^\infty(U).
\]
Here $\widehat\otimes_\pi$ denotes the completion
of the algebraic tensor product with respect to
the projective tensor norm. (Explicitly,
$\|r\|_{\mathfrak A(V,U)}=
\inf_{r=\sum_{n=1}^M\phi_n\otimes\psi_n}
\sum_{n=1}^M\|\phi_n\|_\infty\|\psi_n\|_\infty$
and $\|R\|_{\mathfrak B(V,U;H)}=
\inf_{R=\sum_{n=1}^M\phi_n\otimes T_n\otimes\psi_n}
\sum_{n=1}^M
\|\phi_n\|_\infty
\|T_n\|_{\mathcal L(H)}
\|\psi_n\|_\infty$.)
We use the natural realization of elements of
$\mathfrak A(V,U)$ and $\mathfrak B(V,U;H)$
as essentially bounded measurable scalar-valued and
$\mathcal L(H)$-valued functions on $V\times U$,
respectively, determined by
\[
(\phi\otimes\psi)(v,u)=\phi(v)\psi(u),
\qquad
(\phi\otimes T\otimes\psi)(v,u)
=\phi(v)\psi(u)T.
\]

For $X\in\{U,V\}$ and $\theta\in L^\infty(X)$,
let $M_\theta\in\mathcal L(L^2(X;H))$ be the
multiplication operator
$(M_\theta f)(z)=\theta(z)f(z)$.
For $T\in\mathcal L(H)$, we use the same symbol
for the operator on $L^2(V;H)$ given by
$(Tf)(v)=T(f(v))$.
Given a bounded operator
$S:L^2(U;H)\to L^2(V;H)$ and a finite tensor
$R=\sum_{n=1}^M\phi_n\otimes T_n\otimes\psi_n$,
define
\[
R\odot S
:=
\sum_{n=1}^M M_{\phi_n}T_nSM_{\psi_n}.
\]
By the universal property of the algebraic tensor product,
the trilinear map
$(\phi,T,\psi)\mapsto M_\phi TSM_\psi$
induces a unique linear map $R\mapsto R\odot S$.
Thus the definition is independent of the chosen
finite tensor representation of $R$.

We record two facts for use in the proofs below.

\begin{proposition}\label{prop:projective-algebras}
The spaces $\mathfrak A(V,U)$ and $\mathfrak B(V,U;H)$
are unital Banach algebras under pointwise multiplication.
\end{proposition}

\begin{proof}
The Banach algebra structures follow by applying
the projective tensor product construction
\cite[Lemma~1.5.2 and the following paragraph,
pp.~30--31]{Kaniuth}.
The respective units are $1\otimes1$ and
$1\otimes\Id_H\otimes1$.
\end{proof}

\begin{proposition}\label{prop:projective-schur}
Let $S:L^2(U;H)\to L^2(V;H)$ be bounded.
\begin{enumerate}[label=\textup{(\arabic*)}]
\item
The map $R\mapsto R\odot S$ extends uniquely
by continuity to $\mathfrak B(V,U;H)$.
For every $R\in\mathfrak B(V,U;H)$,
\begin{equation}\label{eq:operator-schur}
\|R\odot S\|_{L^2(U;H)\to L^2(V;H)}
\le
\|R\|_{\mathfrak B(V,U;H)}\|S\|.
\end{equation}

\item
If $S$ has a strongly measurable, essentially bounded
kernel $K_S:V\times U\to\mathcal L(H)$, then
\begin{equation}\label{eq:operator-schur-kernel}
((R\odot S)f)(v)
=
\int_U R(v,u)K_S(v,u)f(u)\,d\mu_U(u)
\end{equation}
for every $R\in\mathfrak B(V,U;H)$,
$f\in L^2(U;H)$, and almost every $v\in V$.
\end{enumerate}
\end{proposition}

\begin{proof}
\textup{(1)}
The assertion follows by iterating the
linearization theorem for bounded bilinear maps
\cite[Theorem~2.9, p.~22]{Ryan}, since
$\|M_\phi TSM_\psi\|
\le
\|\phi\|_\infty\|T\|\|\psi\|_\infty\|S\|$.

\textup{(2)}
Fix $R\in\mathfrak B(V,U;H)$ and $f\in L^2(U;H)$. Since $R$ and $K_S$ are strongly measurable and essentially
bounded, and $L^2(U;H)\subset L^1(U;H)$,
the integral in \eqref{eq:operator-schur-kernel}
is a well-defined Bochner integral for almost every $v\in V$.
Choose finite tensors $R_m$ such that
$\|R_m-R\|_{\mathfrak B(V,U;H)}\to0$.
Expanding the defining finite sum gives
\[
((R_m\odot S)f)(v)
=
\int_U R_m(v,u)K_S(v,u)f(u)\,d\mu_U(u)
\]
for almost every $v\in V$.

Since $U$ and $V$ are probability spaces,
\[
\begin{aligned}
&\left\|
\int_U
\bigl(R_m(\cdot,u)-R(\cdot,u)\bigr)
K_S(\cdot,u)f(u)\,d\mu_U(u)
\right\|_{L^2(V;H)}\\
&\qquad\le
\|R_m-R\|_{\mathfrak B(V,U;H)}
\|K_S\|_\infty
\|f\|_{L^2(U;H)}
\longrightarrow0.
\end{aligned}
\]
By \textup{(1)}, $(R_m\odot S)f\to(R\odot S)f$
in $L^2(V;H)$.
Comparing the two limits gives
\[
((R\odot S)f)(v)
=
\int_U R(v,u)K_S(v,u)f(u)\,d\mu_U(u)
\]
for almost every $v\in V$.

\end{proof}

The following estimate will be used in the proof of
\cref{prop:digit} to bound the projective tensor norms
of operator-valued exponentials.

\begin{lemma}
	\label{lem:schur-exponential}
	Let $U,V$ be probability spaces and $H$ a Banach space.
	For every $h\in\mathfrak A(V,U)$ and
	$T\in\cL(H)$, the operator-valued function
	$(v,u)\longmapsto e^{h(v,u)T}$
	is represented by an element of $\mathfrak B(V,U;H)$,
	and
	$
	\|e^{hT}\|_{\mathfrak B(V,U;H)}
	\le
	\exp\bigl(
	\|h\|_{\mathfrak A(V,U)}
	\|T\|_{\cL(H)}
	\bigr).
	$
\end{lemma}

\begin{proof}
	For fixed $T\in\cL(H)$, the map $\phi\otimes\psi
	\longmapsto
	\phi\otimes T\otimes\psi$
	extends to a bounded linear map $\iota_T:\mathfrak A(V,U)\to\mathfrak B(V,U;H)$
	with $\|\iota_T\|\le\|T\|_{\cL(H)}$.
	Put $X=\iota_T(h)$. Then
	\[
	\|X\|_{\mathfrak B(V,U;H)}
	\le
	\|h\|_{\mathfrak A(V,U)}
	\|T\|_{\cL(H)}.
	\]
	By \cref{prop:projective-algebras}, the exponential series
    defines $e^X\in\mathfrak B(V,U;H)$, and
	\[
	\|e^X\|_{\mathfrak B(V,U;H)}
	\le
	\exp\bigl(\|X\|_{\mathfrak B(V,U;H)}\bigr).
	\]
	Under the natural realization as
	operator-valued functions on $V\times U$,
	the element $e^X$ corresponds to $(v,u)\longmapsto e^{h(v,u)T}$.
	The desired estimate follows.
\end{proof}

On the Schwartz space $\mathcal S(\R)$, we use the Fourier transform
and inverse Fourier transform given by
\[
\widehat f(\xi)=\int_\R e^{-it\xi}f(t)\,dt,
\qquad
\check\psi(t)=\frac{1}{2\pi}\int_\R e^{it\xi}\psi(\xi)\,d\xi.
\]
For a Banach space $E$ and $\Delta>0$, define
\[
\cB_\Delta(E)
=
\bigl\{
f\in L^2(\R;E):
\operatorname{supp}\widehat f
\subseteq[-\Delta,\Delta]
\bigr\},
\]
where the Fourier transform and its support are
understood in the distributional sense.
Let $D=d/dt$ denote distributional differentiation,
and define the translation operators by
$(\tau_hf)(t)=f(t+h)$ for $h\in\R$. To apply \cref{lem:schur-exponential}, we next show that
$D$ is bounded on $\cB_\Delta(E)$ and that translations
on this space are given by $e^{hD}$.

\begin{lemma}\label{lem:derivative}
	The operator $D$ restricts to a bounded linear
	operator on $\cB_\Delta(E)$, satisfying
	$\|D\|_{\cL(\cB_\Delta(E))}\le C\Delta$,
	where $C>0$ is independent of $E$ and $\Delta$.
	Moreover, for every $h\in\R$,
	\[
	\tau_h|_{\cB_\Delta(E)}
	=
	e^{hD}
	=
	\sum_{n=0}^\infty\frac{h^nD^n}{n!},
	\]
	where the series converges in operator norm
	on $\cB_\Delta(E)$.
\end{lemma}

\begin{proof}
	
	Let $X=\cB_\Delta(E)$. Choose
	$\chi\in C_c^\infty(\R)$ equal to one on a neighbourhood
	of $[-1,1]$, and put
	$\kappa_\Delta
	=\mathcal F^{-1}\bigl[\chi(\,\cdot\,/\Delta)\bigr]$.
	Since $\chi(\,\cdot\,/\Delta)$ equals one on a
	neighbourhood of $[-\Delta,\Delta]$, for every $f\in X$
	we have, in the sense of $E$-valued tempered
	distributions, $\widehat{\kappa_\Delta*f}=
	\chi(\,\cdot\,/\Delta)\widehat f=
	\widehat f$
	and $\widehat{\kappa_\Delta'*f}
	=
	i\xi\chi(\xi/\Delta)\widehat f
	=
	i\xi\widehat f
	=
	\widehat{Df}$.
	Since $\kappa_\Delta,\kappa_\Delta'\in L^1(\R)$,
	vector-valued Young's inequality shows that both
	convolutions belong to $L^2(\R;E)$. Consequently,
	$f=\kappa_\Delta*f$ and $Df=\kappa_\Delta'*f$
	in $L^2(\R;E)$.
	Moreover, $\widehat{Df}=i\xi\widehat f$ in the distributional sense. Hence
	$\operatorname{supp}\widehat{Df}
	\subseteq[-\Delta,\Delta]$, so $Df\in X$. Writing $\kappa=\mathcal F^{-1}\chi$, the Fourier scaling formula yields $\kappa_\Delta(t)=\Delta\kappa(\Delta t)$ and 
   $\|\kappa_\Delta'\|_{L^1(\R)}
	=
	\Delta\|\kappa'\|_{L^1(\R)}
	=
	C_\chi\Delta$.
	Consequently, vector-valued Young's inequality gives
	$\|Df\|_{L^2(\R;E)}=\|\kappa_\Delta'*f\|_{L^2(\R;E)}\le
		\|\kappa_\Delta'\|_{L^1(\R)} \|f\|_{L^2(\R;E)}\le
		C_\chi\Delta\|f\|_{L^2(\R;E)}$.
	Thus $D\in\cL(X)$ and $\|D\|_{\cL(X)}\le C_\chi\Delta$.
	
	We now consider translations with the convention
	$(\tau_hf)(t)=f(t+h)$. Since $\widehat{\tau_hf}=e^{ih\xi}\widehat f$
	in the distributional sense, we have $\tau_hf\in X$
	for every $h\in\R$. Furthermore, $\tau_hf=(\tau_h\kappa_\Delta)*f$.
	Since $\kappa_\Delta$ is a Schwartz function, the map
	$h\mapsto\tau_h\kappa_\Delta$ is continuously
	differentiable in $L^1(\R)$, with derivative
	$\tau_h\kappa_\Delta'$. Indeed,
	\[
	\left\|
	\frac{\tau_{h+\varepsilon}\kappa_\Delta
		-\tau_h\kappa_\Delta}{\varepsilon}
	-\tau_h\kappa_\Delta'
	\right\|_{L^1(\R)}
	\le
	\frac{|\varepsilon|}{2}
	\|\kappa_\Delta''\|_{L^1(\R)}
	\longrightarrow0.
	\]
	By vector-valued Young's inequality, it follows that
	$\frac{d}{dh}\tau_hf
	=
	(\tau_h\kappa_\Delta')*f
	=
	\tau_h(Df)
	=
	D(\tau_hf)$ in $L^2(\R;E)$.
	Since $X$ is a closed subspace of $L^2(\R;E)$
	and the difference quotients belong to $X$,
	the derivative also exists in $X$.
	
	Now fix $f\in X$ and define $g(h)=\tau_hf$. Then $g:\R\to X$ is continuously differentiable and satisfies the Banach-space Cauchy problem $g'(h)=Dg(h)$, $g(0)=f$.
	Since $D\in\cL(X)$, the existence and uniqueness
	result for linear Cauchy problems with bounded
	operators
	\cite[Chapter~10, Sections~10.1--10.2]{SchechterFA}
	gives $g(h)=e^{hD}f$ for $h\ge0$. Applying the same result to
	$h\mapsto g(-h)$, with generator $-D$, yields
	the identity for $h\le0$ as well. Therefore $\tau_hf=e^{hD}f$
	for every $f\in X$ and $h\in\R$.
	Finally, since
	$\|D^n\|_{\cL(X)}\le\|D\|_{\cL(X)}^n$,
	we have
	\[\left\|
		\tau_h|_X-
		\sum_{n=0}^N\frac{h^nD^n}{n!}
		\right\|_{\cL(X)}\le
		\sum_{n=N+1}^\infty
		\frac{(|h|C_\chi\Delta)^n}{n!}
		\longrightarrow0
		\qquad(N\to\infty).
	\]
	Thus $\tau_h|_X
	=
	e^{hD}
	=
	\sum_{n=0}^\infty\frac{h^nD^n}{n!}$
	in operator norm.
\end{proof}

\subsection{Reduction to finite Fourier transforms}
Fix integers $A\ge2$ and $N\ge1$, and set $\Lambda_N=2\pi A^N$.
For $\beta\in[0,1]$, define $(\mathcal T_{N,\beta}g)(y,t)
 =\int_0^1e^{i\Lambda_Nxy}g(x,t+\beta xy)\,dx$.
We next bound the operator norm of $\mathcal T_{N,\beta}$
in terms of $\rho_A(E)$.
The proof uses the finite Fourier factorization method
underlying the non-Archimedean estimates in
\cite[Lemma~3.2]{LafforgueBanachT} and
\cite[Proposition~3.1 and Corollary~3.2]{LiaoFFT}.

\begin{proposition}\label{prop:digit}

Let $E$ be a Banach space. Suppose that $\Delta>0$, $\delta\ge0$,
and $\beta\in[0,1]$ satisfy $\beta\Delta\le\delta\Lambda_N$.
Then every $g\in L^2([0,1];\cB_\Delta(E))$ satisfies
\begin{equation}\label{eq:digit-bound}
  \norm{\mathcal T_{N,\beta}g}
  \le C_0\left(C_1e^{C_2\delta A}\rho_A(E)\right)^N\norm g,
\end{equation}
where $C_0,C_1,C_2>0$ are absolute constants and both norms
are taken in $L^2([0,1];\cB_\Delta(E))$.

\end{proposition}

\begin{proof}

Since the endpoints have Lebesgue measure zero,
we may work with $x,y\in[0,1)$. Write
\begin{equation}\label{eq:digits}
	x=\sum_{i=1}^Na_iA^{-i}+A^{-N}u,
	\qquad
	y=\sum_{j=1}^Nb_jA^{-j}+A^{-N}v,
\end{equation}
where $a_i,b_j\in\{0,\ldots,A-1\}$ and
$u,v\in[0,1)$. These are the base-$A$
expansions of $x$ and $y$ up to the first
$N$ places, with $u$ and $v$ representing
the rescaled remainders.

Let
\[
U=(\Z/A\Z)^N_a\times[0,1)_u,
\qquad
V=(\Z/A\Z)^N_b\times[0,1)_v,
\]
be two copies of the same probability space,
equipped with the product of normalized counting
measure and Lebesgue measure. The subscripts
indicate the corresponding coordinates
$(a,u)\in U$ and $(b,v)\in V$, arising from
the parametrizations of $x$ and $y$,
respectively.

For every nonnegative measurable function
$F:[0,1)\to[0,\infty]$, we have
\[
\int_0^1 F(x)\,dx
=
\frac{1}{A^N}
\sum_{a_1,\ldots,a_N=0}^{A-1}
\int_0^1
F\left(
\sum_{i=1}^Na_iA^{-i}+A^{-N}u
\right)\,du,
\]
and similarly,
\[
\int_0^1 F(y)\,dy
=
\frac{1}{A^N}
\sum_{b_1,\ldots,b_N=0}^{A-1}
\int_0^1
F\left(
\sum_{j=1}^Nb_jA^{-j}+A^{-N}v
\right)\,dv.
\]
	
For $j=1,\ldots,N$, put $i_j=N+1-j$,
$m_j=\sum_{k<j}b_kA^{j-1-k}$, and
$\theta_j=\sum_{k>j}b_kA^{j-k}+vA^{j-N}$.	
Then $0\le\theta_j<1$ and
$A^{j-1}y=m_j+\frac{b_j+\theta_j}{A}$.

Consequently,
\[
A^Nxy
=\sum_{j=1}^N a_{i_j}m_j
+\sum_{j=1}^N\frac{a_{i_j}b_j}{A}
+\sum_{j=1}^N\frac{a_{i_j}}A\theta_j+uy.
\]
Write $xy=h_0+h_1$, where $h_0(a,b)=A^{-N}\sum_{j=1}^N a_{i_j}m_j(b)$
and $h_1
=
A^{1-N}\sum_{j=1}^N
\frac{a_{i_j}}A
\frac{b_j+\theta_j}{A}
+
A^{-N}uy$.
The first sum in the expression for $A^Nxy$
is integer-valued.

We isolate the Fourier core.  Let $H=\cB_\Delta(E)$ and define
\[
 \mathcal X_j=L^2\bigl((\Z/A\Z)^j_b\times
 (\Z/A\Z)^{N-j}_a\times[0,1]_u;H\bigr),
 \qquad0\le j\le N.
\]

At stage $j$, the previously introduced output variables
are $b_1,\ldots,b_{j-1}$, and the current input variable
is $a_{i_j}$. Since
$
m_j=\sum_{k=1}^{j-1}b_kA^{j-1-k}
$
depends only on the previously introduced output variables,
we can define the diagonal translation operator
$D_j:\mathcal X_{j-1}\to\mathcal X_{j-1}$ by
\[
\begin{aligned}
	&(D_jf)(b_1,\ldots,b_{j-1},
	a_1,\ldots,a_{i_j},u;t)\\
	&\quad=
	f(b_1,\ldots,b_{j-1},
	a_1,\ldots,a_{i_j},u;
	t+\beta a_{i_j}A^{-N}m_j).
\end{aligned}
\]
For each fixed choice of the discrete variables and $u$,
the shift $\beta a_{i_j}A^{-N}m_j$ is independent of $t$.
Since translations preserve $H=\cB_\Delta(E)$ and are
isometries with respect to its inherited $L^2(\R;E)$ norm,
we obtain
$
\|D_jf\|_{\mathcal X_{j-1}}
=
\|f\|_{\mathcal X_{j-1}}.
$
Thus $D_j$ is an isometry.
Define $F_j:\mathcal X_{j-1}\to\mathcal X_j$ by
\[
\begin{aligned}
	&(F_jf)(b_1,\ldots,b_j,
	a_1,\ldots,a_{i_j-1},u)\\
	&\quad=
	\frac1A\sum_{r=0}^{A-1}
	e^{2\pi i r b_j/A}
	f(b_1,\ldots,b_{j-1},
	a_1,\ldots,a_{i_j-1},r,u),
\end{aligned}
\]
for $f\in\mathcal X_{j-1}$, where $r$ replaces
the input variable $a_{i_j}$. Thus $F_j$
is the normalized finite Fourier transform
from $a_{i_j}$ to $b_j$.
Since $H=\cB_\Delta(E)$ is a closed subspace
of $L^2(\R;E)$, \cref{lem:bochner} gives
$\rho_A(H)\le\rho_A(E)$.
Applying this estimate to the variable
$a_{i_j}$ and integrating over the remaining
variables, we obtain
$
\|F_j\|_{\cL(\mathcal X_{j-1},\mathcal X_j)}
\le \rho_A(H)\le\rho_A(E).
$

Set $Y_N=L^2\bigl((\Z/A\Z)^N_b\times[0,1]_v;H\bigr)$
and define $J:\mathcal X_N\to Y_N$ by $(Jf)(b,v)=\int_0^1 f(b,u)\,du$.
By the Cauchy--Schwarz inequality, we have 
$\|Jf\|_{Y_N}^2
	\le
	\frac1{A^N}\sum_b
	\int_0^1\|f(b,u)\|_H^2\,du
	=
	\|f\|_{\mathcal X_N}^2$. Thus $\|J\|\le1$.
Define the Fourier core operator
\[
S_N=JF_ND_N\cdots F_1D_1:
\mathcal X_0\longrightarrow Y_N.
\]
Write $a=(a_1,\ldots,a_N)$ and
$b=(b_1,\ldots,b_N)$, and recall that
$h_0(a,b)=A^{-N}\sum_{j=1}^N a_{i_j}m_j(b)$.
Expanding the finite Fourier transforms
and using $\tau_s\tau_t=\tau_{s+t}$,
we obtain, for $f\in\mathcal X_0$,
\[
	(S_Nf)(b,v)=
	\frac1{A^N}
	\sum_{a\in(\Z/A\Z)^N}
	\int_0^1
	\exp\left(
	2\pi i\sum_{j=1}^N
	\frac{a_{i_j}b_j}{A}
	\right)
	\tau_{\beta h_0(a,b)}f(a,u)\,du.
\]
Consequently, with respect to the normalized
counting measure on $(\Z/A\Z)^N_a$ and
Lebesgue measure on $[0,1]_u$, the
operator-valued kernel of $S_N$ is
\[
	K((b,v),(a,u))=\exp\left(2\pi i\sum_{j=1}^N
	\frac{a_{i_j}b_j}{A}
	\right)
	\tau_{\beta h_0(a,b)}\in\mathcal L(H).
\]

Finally, since $\|J\|\le1$,
$\|D_j\|=1$, and
$\|F_j\|\le\rho_A(E)$ for every
$j=1,\ldots,N$, the submultiplicativity
of the operator norm gives
\begin{equation}\label{eq:core-bound}
	\|S_N\|
	\le
	\|J\|\prod_{j=1}^N
	\|F_j\|\|D_j\|
	\le \rho_A(E)^N.
\end{equation}

We next estimate the remaining factors in the
integral kernel. Define $\Phi
=
2\pi\left(
\sum_{j=1}^N\frac{a_{i_j}}A\theta_j+uy
\right)$.
Since $a_{i_j}/A$ and $u$ depend only on $(a,u)$,
whereas $\theta_j$ and $y$ depend only on $(b,v)$,
we have the following identity in $\mathfrak A(V,U)$:
$
\Phi
=
2\pi\left(
\sum_{j=1}^N
\theta_j\otimes\frac{a_{i_j}}A
+
y\otimes u
\right)$.
All scalar factors in this representation have
$L^\infty$-norm at most one. Consequently,
\begin{equation}\label{eq:phase-projective}
	\begin{aligned}
		\|\Phi\|_{\mathfrak A(V,U)}
		&\le 2\pi(N+1),\\
		\|e^{i\Phi}\|_{\mathfrak A(V,U)}
		&\le
		\exp\bigl(\|\Phi\|_{\mathfrak A(V,U)}\bigr)
		\le e^{2\pi(N+1)}.
	\end{aligned}
\end{equation}
The second estimate follows from the exponential series
and \cref{prop:projective-algebras}.
Similarly, recalling the definition of $h_1$, we have
$h_1
=
A^{1-N}\sum_{j=1}^N
\frac{b_j+\theta_j}{A}
\otimes\frac{a_{i_j}}A
+
A^{-N}y\otimes u$.
Since $0\le\frac{a_{i_j}}A<1$ and $0\le\frac{b_j+\theta_j}{A}<1$,
we get
\begin{equation}\label{eq:h1-projective}
	\|h_1\|_{\mathfrak A(V,U)}
	\le NA^{1-N}+A^{-N}
	\le (N+1)A^{1-N}.
\end{equation}

By \cref{lem:derivative},
$D=d/dt$ belongs to $\cL(H)$ and satisfies
$\|D\|_{\cL(H)}\le C\Delta$.
Applying \cref{lem:schur-exponential}
with $h=\beta h_1$ and $T=D$, we obtain
\[
\begin{aligned}
	\|e^{\beta h_1D}\|_{\mathfrak B(V,U;H)}
	&\le
	\exp\bigl(
	\beta\|h_1\|_{\mathfrak A(V,U)}
	\|D\|_{\cL(H)}
	\bigr)
	\le
	\exp\bigl(
	C\beta\Delta\|h_1\|_{\mathfrak A(V,U)}
	\bigr).
\end{aligned}
\]
Since $\beta\Delta\le\delta\Lambda_N$
and $\Lambda_N=2\pi A^N$, it follows from
\eqref{eq:h1-projective} that
\begin{equation}\label{eq:translation-tail}
	\begin{aligned}
		\|e^{\beta h_1D}\|_{\mathfrak B(V,U;H)}
		&\le
		\exp\bigl(2\pi C\delta(NA+1)\bigr)
		\le
		\exp\bigl(C'\delta A(N+1)\bigr).
	\end{aligned}
\end{equation}

Let $\mathscr U_N\colon L^2([0,1);H)\to L^2(U;H)$ and
$\mathscr V_N\colon L^2([0,1);H)\to L^2(V;H)$
be the pullback operators associated with
the parametrizations in \eqref{eq:digits}, given by
$(\mathscr U_Ng)(a,u)=g(x(a,u))$ and
$(\mathscr V_Ng)(b,v)=g(y(b,v))$, respectively.
By the change-of-variables formulas above,
both operators are isometric isomorphisms.
By \cref{prop:projective-algebras},
$R:=e^{i\Phi}e^{\beta h_1D}\in\mathfrak B(V,U;H)$,
where $e^{i\Phi}$ is identified with its
natural scalar embedding into $\mathfrak B(V,U;H)$.
Since the kernel of $S_N$ is bounded and takes only
finitely many values, by \cref{prop:projective-schur}\textup{(2)},
for every $f\in L^2(U;H)$ we have
\[
\begin{aligned}
	&((R\odot S_N)f)(b,v)\\
	&\quad=
	\frac1{A^N}\sum_a\int_0^1
	e^{i\Phi}
	e^{2\pi i\sum_j a_{i_j}b_j/A}
	e^{\beta h_1D}\tau_{\beta h_0}
	f(a,u)\,du\\
	&\quad=
	\frac1{A^N}\sum_a\int_0^1
	e^{i\Lambda_Nxy}
	\tau_{\beta xy}f(a,u)\,du\\
	&\quad=
	(\mathscr V_N\mathcal T_{N,\beta}
	\mathscr U_N^{-1}f)(b,v),
\end{aligned}
\]
where $x=x(a,u)$ and $y=y(b,v)$.
Here the second equality follows from $e^{i\Phi}
e^{2\pi i\sum_j a_{i_j}b_j/A}
=e^{i\Lambda_Nxy}$ and $e^{\beta h_1D}\tau_{\beta h_0}=\tau_{\beta(h_0+h_1)}
=\tau_{\beta xy}$.
The last equality follows from the definitions
of the pullback operators and the change of
variables in \eqref{eq:digits}.
Since $f$ is arbitrary, we conclude that
\[
	\mathscr V_N\mathcal T_{N,\beta}\mathscr U_N^{-1}
	=
	\bigl(e^{i\Phi}e^{\beta h_1D}\bigr)\odot S_N
\]
as operators from $L^2(U;H)$ to $L^2(V;H)$.
Combining \eqref{eq:operator-schur},
\eqref{eq:core-bound},
\eqref{eq:phase-projective}, and
\eqref{eq:translation-tail}, we obtain
\[
\begin{aligned}
	\|\mathcal T_{N,\beta}\|_{\cL(L^2([0,1];H))}
	&=
	\|R\odot S_N\|\\
	&\le
	\|e^{i\Phi}\|_{\mathfrak A(V,U)}
	\|e^{\beta h_1D}\|_{\mathfrak B(V,U;H)}
	\|S_N\|\\
	&\le
	e^{2\pi(N+1)}
	e^{C'\delta A(N+1)}
	\rho_A(E)^N\\
	&\le
	\left(
	e^{4\pi}e^{2C'\delta A}\rho_A(E)
	\right)^N.
\end{aligned}
\]
This proves \eqref{eq:digit-bound}.

\end{proof}

\begin{corollary}\label{cor:narrow}
	Let $\mathcal E$ be a class of type $>1$.
	There exist an integer $A\ge2$ and constants
	$\delta>0$, $q\in(0,1)$, and $C>0$ such that
	for every $E\in\mathcal E$, $N\ge1$,
	$\beta\in[0,1]$, and $\Delta>0$ with
	$\beta\Delta\le\delta\Lambda_N$, one has
	\[
	\|\mathcal T_{N,\beta}g\|
	\le Cq^N\|g\|
	\]
	for all $g\in L^2([0,1];\cB_\Delta(E))$.
\end{corollary}

\begin{proof}
	
By \cref{thm:liao}, we can choose $A$ sufficiently large
so that $C_1\sup_{E\in\mathcal E}\rho_A(E)<1/4$.
Then choose $\delta>0$ sufficiently small so that
$e^{C_2\delta A}<2$.
The conclusion follows from \cref{prop:digit}
with $q=1/2$ and $C=C_0$.

\end{proof}

\subsection{Decay estimates for a model integral operator}
We use \cref{cor:narrow} to prove the operator norm
decay estimate in \cref{thm:flat} for kernels defined
by $t-t'+cxy$.
This model arises in \cref{sec:curved}, where suitable
local coordinates reduce the function defining the kernel
to $t-t'+cxy$ plus a remainder.

Fix a compact annulus
$\Omega=\{\xi\in\R:a\le|\xi|\le b\}$,
where $0<a<b<\infty$.
For $\psi\in C_c^\infty(\R)$ with
$\supp\psi\subseteq\Omega$, set $k=\check\psi$ and
$k_M(s)=Mk(Ms)$ for $M\ge1$.

For $c\ne0$, define
\begin{equation}\label{eq:flat-operator}
(\mathsf S_{M,c,\psi}F)(y,t)
=
\int_\R\int_\R
k_M(t-t'+cxy)F(x,t')\,dt'\,dx.
\end{equation}
Equivalently,
\[
(\mathsf S_{M,c,\psi}F)(y,t)
=
\int_\R
\bigl(k_M*F(x,\cdot)\bigr)(t+cxy)\,dx.
\]
Convolution with $k_M$ in the second variable
multiplies the Fourier transform by $\psi(\xi/M)$,
so $k_M*F(x,\cdot)$ has Fourier support contained in
$M\Omega$.

We shall apply this operator on fixed bounded rectangles.
Fix bounded intervals $I_x,I_y,I_t,I_{t'}$ and regard
$\mathsf S_{M,c,\psi}$ as an operator from
$L^2(I_x\times I_{t'})$ to $L^2(I_y\times I_t)$.
The input is extended by zero outside
$I_x\times I_{t'}$, and the output is restricted
to $I_y\times I_t$.

Later, in \cref{cor:flat-amplitude}, we will apply
the estimate below with different choices of $\psi$
and sum the resulting bounds.
We therefore need to specify how the bound depends
on $\psi$.
For each integer $m\ge0$, set
\[
\mathcal N_m(\psi)
=
\sum_{j=0}^m\|\psi^{(j)}\|_\infty.
\]

\begin{theorem}\label{thm:flat}
Let $\mathcal E$ be a class of type $>1$.
With $\Omega$ and the bounded intervals
$I_x,I_y,I_t,I_{t'}$ fixed as above, let
$0<c_0<c_1<\infty$.
There exist $\varepsilon>0$, an integer $m\ge0$
and $0<C<\infty$ depending only on
$\mathcal E$, $\Omega$, the four intervals and $c_0,c_1$
such that for every $E\in\mathcal E$, $M\ge1$,
$c\in\R$ with $c_0\le|c|\le c_1$, and
$\psi\in C_c^\infty(\R)$ with
$\operatorname{supp}\psi\subseteq\Omega$, we have
\[
\bigl\|\mathsf S_{M,c,\psi}\otimes\Id_E\bigr\|_{
\cL\bigl(
L^2(I_x\times I_{t'};E),
L^2(I_y\times I_t;E)
\bigr)}
\le C\mathcal N_m(\psi)M^{-\varepsilon}.
\]
\end{theorem}

\begin{proof}

We first reduce the desired estimate to the case
$x,y\in[0,1]$ and $t,t'\in\R$.
Fix $E\in\mathcal E$, and let
\[
\mathsf S:
L^2(I_x\times\R;E)
\longrightarrow L^2(I_y\times\R;E)
\]
be defined by \eqref{eq:flat-operator}, with
the $x$-integration over $I_x$, the output
restricted to $I_y$, and no restrictions
on the central variables.
Let
\[
\mathsf S^0:
L^2([0,1]\times\R;E)
\longrightarrow L^2([0,1]\times\R;E)
\]
be defined by the same formula with
$x,y\in[0,1]$, and set $L_E=\|\mathsf S^0\|$.

Choose disjoint interval decompositions
\[
I_x=\bigsqcup_{i=1}^{n_x}J_i,
\qquad
I_y=\bigsqcup_{j=1}^{n_y}K_j,
\]
such that there exist $x_i,y_j\in\R$ satisfying
$A_i:=J_i-x_i\subseteq[0,1]$ and
$B_j:=K_j-y_j\subseteq[0,1]$.
For $F\in L^2(I_x\times\R;E)$, define
\[
G_{ij}(u,s)=
\begin{cases}
F(x_i+u,s+cy_ju),&u\in A_i,\\
0,&u\in[0,1]\setminus A_i.
\end{cases}
\]
The substitutions
$x=x_i+u$, $y=y_j+v$, and $t'=s+cy_ju$ give
$
\bigl(\mathsf S(\mathbf1_{J_i}F)\bigr)(y_j+v,t)
=
(\mathsf S^0G_{ij})
\bigl(v,t+cx_i(y_j+v)\bigr)$ for $v\in B_j$.
Here $\mathbf1_{J_i}$ denotes multiplication
by the indicator function of $J_i$ in
the $x$-variable.
Translation invariance in the second variable yields
$
\|G_{ij}\|_{L^2([0,1]\times\R;E)}
=
\|F\|_{L^2(J_i\times\R;E)}.
$
Moreover, since $B_j\subseteq[0,1]$, we have $
\|\mathsf S(\mathbf1_{J_i}F)\|_{L^2(K_j\times\R;E)}
=
\|\mathsf S^0G_{ij}\|_{L^2(B_j\times\R;E)}
\le
\|\mathsf S^0G_{ij}\|_{L^2([0,1]\times\R;E)}
\le
L_E\|G_{ij}\|_{L^2([0,1]\times\R;E)}
=
L_E\|F\|_{L^2(J_i\times\R;E)}$.
Consequently,
\begin{align*}
\|\mathsf SF\|_{L^2(I_y\times\R;E)}^2
&=
\sum_{j=1}^{n_y}
\left\|
\sum_{i=1}^{n_x}
\mathsf S(\mathbf1_{J_i}F)
\right\|_{L^2(K_j\times\R;E)}^2\\
&\le
n_yL_E^2
\left(
\sum_{i=1}^{n_x}
\|F\|_{L^2(J_i\times\R;E)}
\right)^2\\
&\le
n_xn_yL_E^2
\|F\|_{L^2(I_x\times\R;E)}^2.
\end{align*}

Finally, zero extension from $I_{t'}$ to $\R$
and restriction from $\R$ to $I_t$
have operator norms at most one. Therefore,
\[
\bigl\|\mathsf S_{M,c,\psi}\otimes\Id_E\bigr\|_{
\cL\bigl(
L^2(I_x\times I_{t'};E),
L^2(I_y\times I_t;E)
\bigr)}
\le
\sqrt{n_xn_y}\,
\|\mathsf S^0\|_{
\cL\bigl(L^2([0,1]\times\R;E)\bigr)}.
\]
Since $n_x,n_y$ depend only on the fixed intervals,
it suffices to prove the desired estimate for
$\mathsf S^0$.

To apply \cref{cor:narrow}, we first split $\psi$ into
pieces supported in sufficiently short intervals. Let $\delta$ be as in \cref{cor:narrow}, and fix
$0<\theta<\min\{\delta,1/2\}$.
For each $r\in\Omega$, set $U_r=(r-\theta|r|,r+\theta|r|)$.
Since $\Omega\subseteq\bigcup_{r\in\Omega}U_r$
and $\Omega$ is compact, there exist
$r_1,\ldots,r_J\in\Omega$ such that $\Omega\subseteq\bigcup_{\nu=1}^J U_{r_\nu}$.
A smooth partition of unity gives functions
$\chi_\nu\in C_c^\infty(U_{r_\nu})$ satisfying $\sum_{\nu=1}^J\chi_\nu=1$
on $\Omega$.
These points and functions are chosen using only
$\Omega$ and $\theta$, independently of $M$ and $\psi$.
Set $\psi_\nu=\chi_\nu\psi$. Then $\psi=\sum_{\nu=1}^J\psi_\nu$ and
$\operatorname{supp}\psi_\nu\subseteq U_{r_\nu}$.
Define $R_\nu=Mr_\nu$ and $\Delta_\nu=\theta|R_\nu|$.
Since $r_\nu\in\Omega$, we have $aM\le|R_\nu|\le bM$ and $
\Delta_\nu\le\delta|R_\nu|$.
Moreover, we have $\operatorname{supp}\bigl[\psi_\nu(\cdot/M)\bigr]
\subseteq
[R_\nu-\Delta_\nu,R_\nu+\Delta_\nu]$.
Let $\mathsf S_\nu^0$ denote $\mathsf S^0$ with
$\psi$ replaced by $\psi_\nu$.
By linearity, we get $\mathsf S^0=\sum_{\nu=1}^J\mathsf S_\nu^0$.
For each integer $m\ge0$, the Leibniz rule gives
$
\sum_{\nu=1}^J\mathcal N_m(\psi_\nu)
\le C_m\mathcal N_m(\psi).
$
Here $C_m$ depends only on $\Omega$, $\delta$, and $m$.

Fix $\nu$ and write
$R=R_\nu$, $\Delta=\Delta_\nu$, and $\psi_R=\psi_\nu$.
Define
\[
(P_RF)(x,t)
=
\int_\R
M\check\psi_R\bigl(M(t-t')\bigr)F(x,t')\,dt'.
\]
Taking the Fourier transform in the second variable,
understood in the distributional sense, gives
$\widehat{P_RF}(x,\xi)
=
\psi_R(\xi/M)\widehat F(x,\xi)$.
In particular, $\operatorname{supp}_\xi\widehat{P_RF}
\subseteq[R-\Delta,R+\Delta]$.
Now set $g(x,t)=e^{-iRt}(P_RF)(x,t)$.
Since $\widehat g(x,\xi)=\widehat{P_RF}(x,\xi+R)$, we have
$\operatorname{supp}\widehat g(x,\cdot)\subseteq[-\Delta,\Delta]$
for almost every $x$. By Young's inequality,
$g\in L^2([0,1]\times\R;E)$, so $g\in L^2([0,1];\cB_\Delta(E))$.
By the definition of $P_R$, $(\mathsf S_\nu^0F)(y,t)
=
\int_0^1(P_RF)(x,t+cxy)\,dx$.
Substituting $(P_RF)(x,s)=e^{iRs}g(x,s)$ yields
\[
e^{-iRt}(\mathsf S_\nu^0F)(y,t)
=
\int_0^1e^{icRxy}g(x,t+cxy)\,dx.
\]

For $R,c\ne0$, define
\[
(\mathsf A_{R,c}h)(y,t)
=
\int_0^1 e^{icRxy}h(x,t+cxy)\,dx.
\]
The preceding identity gives
$\|\mathsf S_\nu^0F\|_2=\|\mathsf A_{R,c}g\|_2$.

Let $A$ and $q$ be as in \cref{cor:narrow}, and set
$\varepsilon=-\log q/\log A>0$.
We will prove the following estimate
\begin{equation}\label{eq:flat-auxiliary}
\|\mathsf A_{R,c}h\|_2
\le CM^{-\varepsilon}\|h\|_2
\end{equation}
for every $E\in\mathcal E$, $M\ge1$, $\Delta>0$,
and $h\in L^2([0,1];\cB_\Delta(E))$ satisfying
$aM\le|R|\le bM$, $\Delta\le\delta|R|$,
and $c_0\le|c|\le c_1$.
Here $C$ depends only on
$\mathcal E$, $a,b$, and $c_0,c_1$.

We first explain how this estimate yields the desired bound.
By Young's inequality,
$\|g\|_2=\|P_RF\|_2
\le\|M\check\psi_R(M\cdot)\|_{L^1(\R)}\|F\|_2$.
A change of variables gives
$\|M\check\psi_R(M\cdot)\|_{L^1(\R)}
=\|\check\psi_R\|_{L^1(\R)}$.
To estimate the kernel norm, integration by parts gives
\[
(1+s^2)\check\psi_R(s)
=
\frac1{2\pi}\int_\R
e^{is\xi}\bigl(\psi_R(\xi)-\psi_R''(\xi)\bigr)\,d\xi.
\]
Since $\operatorname{supp}\psi_R\subseteq\Omega$, it follows that
\[
|\check\psi_R(s)|
\le
\frac{\|\psi_R\|_{L^1(\R)}
      +\|\psi_R''\|_{L^1(\R)}}{2\pi(1+s^2)}
\le
\frac{C_\Omega\mathcal N_2(\psi_R)}{1+s^2}.
\]
Integrating in $s$, we obtain
$\|\check\psi_R\|_{L^1(\R)}
\le C_\Omega\mathcal N_2(\psi_R)$.

Assuming \eqref{eq:flat-auxiliary} for the moment,
we therefore have
\[
\|\mathsf S_\nu^0F\|_2
=\|\mathsf A_{R,c}g\|_2
\le CM^{-\varepsilon}\|g\|_2
\le
C\mathcal N_2(\psi_\nu)M^{-\varepsilon}\|F\|_2.
\]
Summing over $\nu$ and using the seminorm estimate
for the decomposition established above gives
\[
\|\mathsf S^0\|
\le\sum_{\nu=1}^J\|\mathsf S_\nu^0\|
\le
CM^{-\varepsilon}
\sum_{\nu=1}^J\mathcal N_2(\psi_\nu)
\le C\mathcal N_2(\psi)M^{-\varepsilon}.
\]
Combining this with the reduction to the unit intervals yields
\[
\bigl\|\mathsf S_{M,c,\psi}\otimes\Id_E\bigr\|_{
\cL\bigl(
L^2(I_x\times I_{t'};E),
L^2(I_y\times I_t;E)
\bigr)}
\le \sqrt{n_xn_y}\,\|\mathsf S^0\|
\le C\mathcal N_2(\psi)M^{-\varepsilon}.
\]
Thus \eqref{eq:flat-auxiliary} implies the theorem with $m=2$.

It remains to prove \eqref{eq:flat-auxiliary}.
For a given $h$, write
$h^{\mathrm t}(x,s)=h(x,-s)$ and
$h^{\mathrm x}(x,s)=h(1-x,s)$.
Direct changes of variables give
\[
\begin{aligned}
(\mathsf A_{-r,c}h)(y,-t)
&=
(\mathsf A_{r,-c}h^{\mathrm t})(y,t),\\
(\mathsf A_{r,-d}h)(y,t)
&=
e^{-idry}
(\mathsf A_{r,d}h^{\mathrm x})(y,t-dy),
\end{aligned}
\]
where $r,d>0$.
We have
$\|h^{\mathrm t}\|_2=\|h^{\mathrm x}\|_2=\|h\|_2$.
Since
$\widehat{h^{\mathrm t}}(x,\xi)=\widehat h(x,-\xi)$
and
$\widehat{h^{\mathrm x}}(x,\xi)=\widehat h(1-x,\xi)$,
both $h^{\mathrm t}$ and $h^{\mathrm x}$
belong to $L^2([0,1];\cB_\Delta(E))$.
Changing variables in $t$ and using $|e^{-idry}|=1$,
we obtain from the two formulas for $\mathsf A$ above
\begin{equation}\label{eq:flat-sign-reduction}
\|\mathsf A_{-r,c}h\|_2
=
\|\mathsf A_{r,-c}h^{\mathrm t}\|_2,
\qquad
\|\mathsf A_{r,-d}h\|_2
=
\|\mathsf A_{r,d}h^{\mathrm x}\|_2.
\end{equation}
Thus it suffices to consider $R=r>0$ and $c=d>0$.

By the triangle and Cauchy--Schwarz inequalities,
followed by Tonelli's theorem and translation
invariance,
\[
\|\mathsf A_{r,d}h\|_2^2
\le
\int_0^1\int_0^1\int_\R
\|h(x,t+dxy)\|_E^2\,dt\,dx\,dy
=\|h\|_2^2.
\]
If $r<2\pi A$, then $aM\le r$ implies
$M<2\pi A/a$, and hence
$\|\mathsf A_{r,d}h\|_2
\le\|h\|_2
\le(2\pi A/a)^\varepsilon M^{-\varepsilon}\|h\|_2$.

Now assume $r\ge2\pi A$.
Choose $N\ge1$ such that
$2\pi A^N\le r<2\pi A^{N+1}$,
and set $\beta=2\pi A^N/r$ and $\lambda=d/\beta$.
Then $A^{-1}<\beta\le1$, $c_0\le\lambda<Ac_1$,
and $\beta\Delta\le\beta\delta r=\delta\Lambda_N$.
The substitution $u=\lambda x$ gives
\[
(\mathsf A_{r,d}h)(y,t)
=
\frac1\lambda
\int_0^\lambda
e^{i\Lambda_Nuy}
h\left(\frac u\lambda,t+\beta uy\right)\,du.
\]

To apply \cref{cor:narrow}, split $[0,\lambda]$
into intervals of length at most one.
Set $K=\lceil\lambda\rceil$ and, for $1\le j\le K$,
define
$B_j=([j-1,j]\cap[0,\lambda])-(j-1)\subseteq[0,1]$
and
\[
h_j(u,s)=
\begin{cases}
h\left(\dfrac{j-1+u}{\lambda},s\right),
&u\in B_j,\\
0,&u\in[0,1]\setminus B_j.
\end{cases}
\]
These functions satisfy
$h_j\in L^2([0,1];\cB_\Delta(E))$ and
$\sum_{j=1}^K\|h_j\|_2^2=\lambda\|h\|_2^2$.
Splitting the integral above and substituting
$u=j-1+v$ in its $j$th part gives
\[
(\mathsf A_{r,d}h)(y,t)
=
\frac1\lambda
\sum_{j=1}^K
e^{i\Lambda_N(j-1)y}
(\mathcal T_{N,\beta}h_j)
\bigl(y,t+\beta(j-1)y\bigr).
\]
By \cref{cor:narrow},
$\|\mathcal T_{N,\beta}h_j\|_2
\le Cq^N\|h_j\|_2$.
Since the factors $e^{i\Lambda_N(j-1)y}$
have modulus one and translations in $t$
preserve the norm, we obtain
\[
\|\mathsf A_{r,d}h\|_2
\le
\frac{Cq^N}{\lambda}
\sum_{j=1}^K\|h_j\|_2
\le
\frac{Cq^N\sqrt K}{\lambda}
\left(\sum_{j=1}^K\|h_j\|_2^2\right)^{1/2}
=
Cq^N\sqrt{\frac K\lambda}\,\|h\|_2
\le C'q^N\|h\|_2,
\]
where the last inequality follows from
$K\le\lambda+1$ and $\lambda\ge c_0$.

Finally, the choice of $N$ and the inequality $r\ge aM$ give
\[
q^N
=
(A^N)^{-\varepsilon}
\le
\left(\frac{2\pi A}{r}\right)^\varepsilon
\le
\left(\frac{2\pi A}{a}\right)^\varepsilon
M^{-\varepsilon}.
\]
This proves \eqref{eq:flat-auxiliary} for positive $R,c$.
By \eqref{eq:flat-sign-reduction} the same estimate
holds for all $R,c\ne0$. This completes the proof.

\end{proof}

\section{Integral operator estimates under a curvature condition}
\label{sec:curved}

\subsection{Integral operators with smooth amplitudes}

Since $k=\check\psi$, Fourier inversion gives
$Mk(Ms)=
\frac M{2\pi}\int_\R
e^{iM\omega s}\psi(\omega)\,d\omega$.
Thus the operator in \cref{thm:flat} can also be written as
\begin{equation}\label{eq:flat-fourier}
(\mathsf S_{M,c,\psi}F)(y,t)
=
\frac M{2\pi}\iiint
e^{iM\omega(t-t'+cxy)}
\psi(\omega)F(x,t')\,d\omega\,dt'\,dx.
\end{equation}
Here the factor $\psi(\omega)$ depends only on $\omega$.
We now consider operators of the form
\begin{equation}\label{eq:flat-amplitude}
(\mathsf A_{M,c,B}F)(y,t)
=
\frac M{2\pi}\iiint_{\R^3}
e^{iM\omega(t-t'+cxy)}
B(y,t,x,t',\omega)F(x,t')
\,d\omega\,dt'\,dx,
\end{equation}
where $B\in C_c^{r}(\R^5)$ satisfies
$\operatorname{supp}B
\subseteq
(I_y\times I_t)\times
(I_x\times I_{t'})\times\Omega$.
The integer ${r}\ge0$ will be chosen in
\cref{cor:flat-amplitude}.
Our goal is to prove
\[
\|\mathsf A_{M,c,B}\otimes\Id_E\|
\le
C\|B\|_{C^r(\R^5)}M^{-\varepsilon},
\]
with the same exponent $\varepsilon$ as in
\cref{thm:flat}. These operators arise from the rescaled expression
\eqref{eq:scaled-block} in the proof of \cref{thm:curved}.

Write $q=(y,t)$ and $p=(x,t')$. If
$B(q,p,\omega)=u(q)v(p)\phi(\omega)$ with bounded $u,v$
and $\phi\in C_c^\infty(\R)$ supported in $\Omega$, then
$\mathsf A_{M,c,B}F=u\,\mathsf S_{M,c,\phi}(vF)$.
Its norm is controlled by \cref{thm:flat} and the bounds
for multiplication by $u$ and $v$.

For a general amplitude $B$, \cref{lem:smooth-projective}
gives a decomposition into products
$u_\nu(q)v_\nu(p)\phi_\nu(\omega)$ with a bound that allows
the resulting operator estimates to be summed.
The proof uses a Fourier-series argument as in
\cite[Lecture~15, pp.~78--79]{RoeAnalysis}.
For the corresponding tensor-product identity for smooth
function spaces, see \cite[Theorem~51.6]{Treves}.

\begin{lemma}
\label{lem:smooth-projective}
Let $Q\subset\R^{d_1}$ and $P\subset\R^{d_2}$
be fixed compact boxes with nonempty interiors,
and let $I\subset\R$ be a fixed compact interval
with nonempty interior.
For every integer $m_0\ge0$, there exist an integer
$m\ge m_0$ and a constant $C>0$, depending only on
$Q,P,I$ and $m_0$, such that every
$B\in C^m(Q\times P\times I)$ admits a representation
\[
B(q,p,\omega)
=
\sum_{\nu=1}^\infty
u_\nu(q)v_\nu(p)\phi_\nu(\omega),
\]
where $u_\nu\in C(Q)$, $v_\nu\in C(P)$,
and $\phi_\nu\in C^{m_0}(I)$ satisfy
\[
\sum_{\nu=1}^\infty
\|u_\nu\|_\infty
\|v_\nu\|_\infty
\|\phi_\nu\|_{C^{m_0}(I)}
\le C\|B\|_{C^m(Q\times P\times I)}.
\]
\end{lemma}

\begin{proof}
Set $D=d_1+d_2+1$ and $m=m_0+D+1$,
and write $K=Q\times P\times I$.
We will obtain the desired products from a Fourier series.
For this purpose, we first construct a periodic $C^m$
function representing $B$, with a controlled $C^m$ norm.

We enlarge the domain so that there is room for a cutoff
which equals one on $K$ and vanishes near the boundary
of the larger box.
Choose fixed compact boxes $Q'$ and $P'$
and a compact interval $I'$ such that
$Q\subset\operatorname{int}Q'$,
$P\subset\operatorname{int}P'$, and
$I\subset\operatorname{int}I'$,
and set $K'=Q'\times P'\times I'$.
Applying the reflection construction in
\cite[Theorem~5.19 and its proof]{AdamsFournier}
successively at the endpoints of the coordinate intervals,
we obtain an extension $\widetilde B\in C^m(K')$
satisfying $\widetilde B|_K=B$ and
\[
\|\widetilde B\|_{C^m(K')}
\le C\|B\|_{C^m(K)}.
\]

To make the periodic copies join smoothly, we arrange
that the extension vanishes near the boundary of $K'$.
Choose a fixed $\chi\in C_c^\infty(\operatorname{int}K')$
with $\chi=1$ on $K$, and set $H=\chi\widetilde B$.
Then $H|_K=B$, $H$ vanishes near the boundary of $K'$,
and the Leibniz rule gives
\[
\|H\|_{C^m(K')}
\le C\|\widetilde B\|_{C^m(K')}
\le C\|B\|_{C^m(K)}.
\]

We now map $K'$ onto the standard periodic cell
$[-\pi,\pi]^D$.
Choose fixed coordinatewise affine bijections
$T_Q:Q'\longrightarrow[-\pi,\pi]^{d_1}$
and $T_P:P'\longrightarrow[-\pi,\pi]^{d_2}$,
and an affine bijection
$s(\omega)=\kappa\omega+\sigma$,
with $\kappa>0$, from $I'$ onto $[-\pi,\pi]$.
On $[-\pi,\pi]^D$, define
\[
F(\xi,\eta,\tau)
=
H\bigl(T_Q^{-1}(\xi),T_P^{-1}(\eta),s^{-1}(\tau)\bigr).
\]
Since $H$ vanishes near the boundary of $K'$,
the values and derivatives of $F$ agree across
opposite faces of this cell.
Its $2\pi$-periodic extension therefore belongs to
$C^m(\mathbb T^D)$.
The coordinate changes are fixed, so the chain rule gives
$\|F\|_{C^m(\mathbb T^D)}\le C\|B\|_{C^m(K)}$.
Since $H|_K=B$, we also have
$B(q,p,\omega)=F(T_Q(q),T_P(p),s(\omega))$
for $(q,p,\omega)\in K$.

It remains to estimate the Fourier coefficients of $F$
with the weights needed to control $m_0$ derivatives
in the frequency variable.
For $n=(\alpha,\beta,k)\in
\Z^{d_1}\times\Z^{d_2}\times\Z$, let
\[
c_{\alpha,\beta,k}
=
\frac1{(2\pi)^D}
\int_{\mathbb T^D}F(z)e^{-in\cdot z}\,dz.
\]
If $n\ne0$, choose a coordinate $j$ such that
$|n_j|=|n|_\infty$.
Integration by parts $m$ times in that coordinate gives
$|c_{\alpha,\beta,k}|
\le |n|_\infty^{-m}\|\partial_j^mF\|_\infty$.
For $n=0$, the coefficient is bounded by $\|F\|_\infty$.
Combining the two cases yields
\[
|c_{\alpha,\beta,k}|
\le
C(1+|n|_\infty)^{-m}
\|F\|_{C^m(\mathbb T^D)}.
\]

Differentiating $e^{ik\tau}$ $j$ times produces
the factor $(ik)^j$.
Thus the weight $(1+|k|)^{m_0}$ controls all derivatives
of this factor through order $m_0$.
Since $|k|\le|n|_\infty$ and $m-m_0>D$, we obtain
\[
\begin{aligned}
\sum_{\alpha,\beta,k}
|c_{\alpha,\beta,k}|(1+|k|)^{m_0}
&\le
C\|F\|_{C^m(\mathbb T^D)}
\sum_{n\in\Z^D}
(1+|n|_\infty)^{-(m-m_0)}\\
&\le C\|B\|_{C^m(K)}.
\end{aligned}
\]
In particular, the Fourier series converges absolutely
and uniformly.
Its sum has the same Fourier coefficients as $F$,
so it equals $F$ by uniqueness. Hence
\[
F(\xi,\eta,\tau)
=
\sum_{\alpha,\beta,k}
c_{\alpha,\beta,k}
e^{i\alpha\cdot\xi}
e^{i\beta\cdot\eta}
e^{ik\tau}.
\]

We now return to the original coordinates.
Since $T_Q$, $T_P$, and $s$ act separately on
$q$, $p$, and $\omega$, each Fourier term remains
a product of functions of these three groups of variables.
For each triple $\nu=(\alpha,\beta,k)$, define
\begin{equation}\label{eq:projective-factors}
u_\nu(q)=e^{i\alpha\cdot T_Q(q)},
\qquad
v_\nu(p)=e^{i\beta\cdot T_P(p)},
\qquad
\phi_\nu(\omega)
=
c_{\alpha,\beta,k}e^{ik s(\omega)}.
\end{equation}
Substituting the Fourier expansion into the identity
relating $B$ to $F$ gives
\[
B(q,p,\omega)
=
\sum_\nu u_\nu(q)v_\nu(p)\phi_\nu(\omega).
\]
The functions $u_\nu$ and $v_\nu$ are restrictions
of smooth exponentials, so they belong to $C(Q)$
and $C(P)$, respectively, and satisfy
$\|u_\nu\|_\infty=\|v_\nu\|_\infty=1$.

It therefore remains to control the frequency factors.
For $0\le j\le m_0$, the identity
$s(\omega)=\kappa\omega+\sigma$ gives
$\phi_\nu^{(j)}(\omega)
=c_{\alpha,\beta,k}(i\kappa k)^j e^{ik s(\omega)}$.
Consequently,
\[
\|\phi_\nu\|_{C^{m_0}(I)}
\le
|c_{\alpha,\beta,k}|
\sum_{j=0}^{m_0}(\kappa|k|)^j
\le
\max\{1,\kappa\}^{m_0}
|c_{\alpha,\beta,k}|(1+|k|)^{m_0}.
\]
The weighted summability estimate established above
now gives
\[
\begin{aligned}
\sum_\nu
\|u_\nu\|_\infty
\|v_\nu\|_\infty
\|\phi_\nu\|_{C^{m_0}(I)}
&\le
C\sum_{\alpha,\beta,k}
|c_{\alpha,\beta,k}|(1+|k|)^{m_0}\\
&\le C\|B\|_{C^m(K)}.
\end{aligned}
\]
All choices and constants depend only on $Q,P,I$
and $m_0$.
Enumerating the triples $(\alpha,\beta,k)$
by $\nu=1,2,\ldots$ gives the required representation.
\end{proof}

We now apply \cref{lem:smooth-projective}
to the operator in \eqref{eq:flat-amplitude}.
Recall the class $\mathcal E$ of type $>1$,
the fixed bounded intervals $I_x,I_y,I_t,I_{t'}$,
and the constants $0<c_0<c_1<\infty$
from \cref{thm:flat}.
We retain the annulus
$\Omega=\{\omega\in\R:a\le|\omega|\le b\}$,
where $0<a<b<\infty$, and the exponent
$\varepsilon>0$ chosen in the proof of that theorem.
Set
\[
D_y=\overline{I_y}\times\overline{I_t},
\qquad
D_x=\overline{I_x}\times\overline{I_{t'}}.
\]
Since the endpoints have measure zero, replacing
the intervals by their closures does not change
the corresponding $L^2$ spaces.

\begin{corollary}
\label{cor:flat-amplitude}
There exist an integer $r\ge0$ and a constant $C>0$,
depending only on $\mathcal E$, $\Omega$, $D_y,D_x$,
and $c_0,c_1$, such that for every
$B\in C_c^r(\R^5)$ with
$\operatorname{supp}B\subseteq D_y\times D_x\times\Omega$,
the operator defined in \eqref{eq:flat-amplitude}
satisfies
\[
\bigl\|\mathsf A_{M,c,B}\otimes\Id_E\bigr\|_{
\cL\bigl(
L^2(D_x;E),
L^2(D_y;E)
\bigr)}
\le
C\|B\|_{C^r(\R^5)}M^{-\varepsilon}
\]
uniformly for $E\in\mathcal E$, $M\ge1$,
and $c\in\R$ with $c_0\le|c|\le c_1$.
\end{corollary}

\begin{proof}We will decompose $B$ into terms of the form
$u_\nu(q)v_\nu(p)\psi_\nu(\omega)$
and apply \cref{thm:flat} to the corresponding
operators.
Choose a fixed
$\chi\in C_c^\infty(\R\setminus\{0\})$
equal to one on $\Omega$, and fix a compact annulus
$\Omega_*$ containing $\operatorname{supp}\chi$.
In the proof of \cref{thm:flat}, the exponent
$\varepsilon=-\log q/\log A$ depends only on
$\mathcal E$.
We may therefore use the same exponent
when applying that theorem with $\Omega_*$.

Let $m_0$ be the derivative order required by
\cref{thm:flat} for $\Omega_*$, and set $r=m_0+6$.
Choose a fixed compact interval $J$ such that
$\operatorname{supp}\chi\subset\operatorname{int}J$.
Apply \cref{lem:smooth-projective} with
$Q=D_y$, $P=D_x$, $d_1=d_2=2$, and $I=J$,
choosing the factors as in \eqref{eq:projective-factors}.
We obtain $B(q,p,\omega)=
\sum_\nu u_\nu(q)v_\nu(p)\phi_\nu(\omega)$,
where each $\phi_\nu$ is smooth and
\[
\sum_\nu
\|u_\nu\|_\infty
\|v_\nu\|_\infty
\|\phi_\nu\|_{C^{m_0}(J)}
\le C\|B\|_{C^r(\R^5)}.
\]
Set $\psi_\nu=\chi\phi_\nu$ on $J$ and extend it
by zero outside $J$.
Since $\chi$ vanishes near the endpoints of $J$
and $\phi_\nu$ is smooth, this extension belongs to
$C_c^\infty(\R)$.
Moreover,
$\operatorname{supp}\psi_\nu
\subseteq\operatorname{supp}\chi\subseteq\Omega_*$.
The Leibniz rule gives
\[
\mathcal N_{m_0}(\psi_\nu)
\le C_{\chi,m_0}\|\phi_\nu\|_{C^{m_0}(J)},
\]
and hence the preceding summability estimate implies
\[
\sum_\nu
\|u_\nu\|_\infty
\|v_\nu\|_\infty
\mathcal N_{m_0}(\psi_\nu)
\le C\|B\|_{C^r(\R^5)}.
\]

Since $B=\chi B$, we have
$B(q,p,\omega)=\sum_\nu u_\nu(q)v_\nu(p)\psi_\nu(\omega)$ with
$(q,p,\omega)\in D_y\times D_x\times\R$.
This series converges absolutely and uniformly,
because
$\|\psi_\nu\|_\infty\le\mathcal N_{m_0}(\psi_\nu)$.

We claim that, for $F\in L^2(D_x;E)$,
\[
(\mathsf A_{M,c,B}\otimes\Id_E)F
=
\sum_\nu
u_\nu\,
(\mathsf S_{M,c,\psi_\nu}\otimes\Id_E)(v_\nu F),
\]
with convergence in operator norm.
To justify this identity, set
$B_N(q,p,\omega)=
\sum_{\nu=1}^N
u_\nu(q)v_\nu(p)\psi_\nu(\omega)$.
By \eqref{eq:flat-fourier} and
\eqref{eq:flat-amplitude} we have
\[
(\mathsf A_{M,c,B_N}\otimes\Id_E)F
=
\sum_{\nu=1}^N
u_\nu\,
(\mathsf S_{M,c,\psi_\nu}\otimes\Id_E)(v_\nu F).
\]
Denote the operator in the $\nu$-th summand
by $T_\nu$.

To identify the sum of the operators, we first
pass to the limit for each fixed $M$.
The decay in $M$ will then follow from
\cref{thm:flat}. Since $B_N\to B$ uniformly and
$B_N(q,p,\omega)=B(q,p,\omega)=0$
whenever $\omega\notin\Omega_*$,
\eqref{eq:flat-amplitude} and the Cauchy--Schwarz
inequality on $D_x$ give
\[
\left\|
\mathsf A_{M,c,B}\otimes\Id_E
-\sum_{\nu=1}^N T_\nu
\right\|
\le
CM\|B-B_N\|_\infty
\longrightarrow0
\]
 where $C$ depends only on
$D_x,D_y$, and $\Omega_*$.
This proves the claimed operator identity.

Finally, \cref{thm:flat} and the bounds for
multiplication operators give
\[
\sum_\nu\|T_\nu\|
\le
CM^{-\varepsilon}
\sum_\nu
\|u_\nu\|_\infty
\|v_\nu\|_\infty
\mathcal N_{m_0}(\psi_\nu)
\le
C\|B\|_{C^r(\R^5)}M^{-\varepsilon}.
\]
Consequently, $\|\mathsf A_{M,c,B}\otimes\Id_E\|
\le
C\|B\|_{C^r(\R^5)}M^{-\varepsilon}$ as required.
\end{proof}

\subsection{Local coordinates and a partition of unity}
The estimate in \cref{cor:flat-amplitude} applies
to the phase $t-t'+cxy$.
We now prepare to use it for kernels containing
$Rk(R\rho(q,p))$, as in \cref{thm:curved}.
Here $R\ge1$, and $Q,P$ are two-dimensional smooth
manifolds equipped with smooth positive densities
$d\mu_Q$ and $d\mu_P$, respectively.
Let $\rho$ be a real-valued function on $Q\times P$.
We seek local coordinates near $\{\rho=0\}$
that relate $\rho$ to the phase $t-t'+cxy$.

To obtain a nonzero $xy$ coefficient in the local
normal form at $z_0\in\{\rho=0\}$,
we require $\Rot(\rho)(z_0)\ne0$.
The rotational curvature is defined in local coordinates by
\begin{equation}\label{eq:rot}
\Rot(\rho)=
\det\begin{pmatrix}
0&(\partial_q\rho)^T\\
\partial_p\rho&\partial^2_{pq}\rho
\end{pmatrix}.
\end{equation}
The condition $\Rot(\rho)\ne0$ on $\{\rho=0\}$
is invariant under changes of product coordinates.
See \cite[Section~1, pp.~109--112]{PhongSteinRadonI}
for the rotational curvature condition.

For the next lemma, fix a compact set $K\subset Q\times P$
and a relatively compact open neighborhood $U$ of $K$.
Choose finitely many coordinate parametrizations
$\sigma_\ell:\Omega_\ell^Q\subset\mathbb R^2\to Q$ and
$\tau_\ell:\Omega_\ell^P\subset\mathbb R^2\to P$,
$\ell=1,\ldots,J$, whose product images cover $\overline U$.
Let $m\ge3$ be an integer and assume that
$\rho|_U\in C^m(U;\R)$.
Set
\[
\widetilde U_\ell
=(\sigma_\ell\times\tau_\ell)^{-1}(U),
\qquad
\widetilde K_\ell
=(\sigma_\ell\times\tau_\ell)^{-1}(K).
\]
On $\widetilde U_\ell$, define
\[
f_\ell(\bar y,\bar t,\bar x,\bar t')
=\rho\bigl(\sigma_\ell(\bar y,\bar t),
\tau_\ell(\bar x,\bar t')\bigr).
\]
Assume that, for fixed $B_0,\kappa>0$
and every $\ell=1,\ldots,J$,
\[
\begin{gathered}
\|f_\ell\|_{C^m(\widetilde U_\ell)}\le B_0,\\
\min\bigl\{
|D_1f_\ell(w)|,\,
|D_2f_\ell(w)|,\,
|\Rot(f_\ell)(w)|
\bigr\}\ge\kappa
\quad
(w\in\widetilde K_\ell\cap\{f_\ell=0\}),
\end{gathered}
\]
where $D_1,D_2$ denote derivatives with respect to
$(\bar y,\bar t)$ and $(\bar x,\bar t')$, respectively.

Under these assumptions, the following lemma constructs
local coordinates in which $\rho$ takes the form
$t-t'+c_{z_0}xy+L\Theta_{z_0,L}$ after the scaling
$(y,t,x,t')\mapsto(Ly,L^2t,Lx,L^2t')$ and division by $L^2$.
Under this scaling, the terms $t$, $t'$, and $xy$
all acquire the same factor $L^2$.

\begin{lemma}\label{lem:normal-form}
Fix a compact box $D\subset\R^4$.
There exist constants $r_0,L_0>0$ and $C\ge1$,
depending only on $K,U,D$, the fixed charts,
$m,B_0$, and $\kappa$, such that the following holds.

For every $z_0=(q_0,p_0)\in K\cap\{\rho=0\}$,
there are $C^m$ coordinate maps
$\mathcal Q_{z_0}:(-r_0,r_0)^2\to Q$ and
$\mathcal P_{z_0}:(-r_0,r_0)^2\to P$
satisfying:
\begin{enumerate}[label=\textup{(\arabic*)}]
\item
$\mathcal Q_{z_0}(0)=q_0$,
$\mathcal P_{z_0}(0)=p_0$, and their product image
is contained in $U$.

\item
For every $0<L\le L_0$ and $(y,t,x,t')\in D$,
\begin{equation}\label{eq:scaled-normal}
L^{-2}\rho\bigl(
\mathcal Q_{z_0}(Ly,L^2t),
\mathcal P_{z_0}(Lx,L^2t')
\bigr)
=t-t'+c_{z_0}yx+L\Theta_{z_0,L}(y,t,x,t'),
\end{equation}
where $C^{-1}\le|c_{z_0}|\le C$ and
$\|\Theta_{z_0,L}\|_{C^{m-3}(D)}\le C$.
\end{enumerate}
\end{lemma}

\begin{proof}
Choose finitely many product coordinate boxes covering $K$,
with closures contained in $U$ and in the fixed
charts. All constants below depend only on the data
specified in the statement.

Fix $z_0=(q_0,p_0)\in K\cap\{\rho=0\}$.
In one of these boxes, write
$q=\sigma(\eta,s)$ and $p=\tau(\xi,s')$,
where
$\sigma:V_Q\subset\R^2\longrightarrow Q$ and
$\tau:V_P\subset\R^2\longrightarrow P$ are the inverses of the two coordinate charts
used to define this product box.
Translate the coordinates so that
$\sigma(0,0)=q_0$ and $\tau(0,0)=p_0$.

After independently interchanging $\eta$ with $s$
and $\xi$ with $s'$ if necessary, the maps
\[
A(\eta,s)=\bigl(\eta,\rho(\sigma(\eta,s),p_0)\bigr),\qquad
B(\xi,s')=\bigl(\xi,-\rho(q_0,\tau(\xi,s'))\bigr)
\]
satisfy $|\det DA(0,0)|\ge\kappa/\sqrt2$ and
$|\det DB(0,0)|\ge\kappa/\sqrt2$.
Indeed, their Jacobian determinants are
\[
\det DA(\eta,s)=
\frac{\partial}{\partial s}
\bigl[\rho(\sigma(\eta,s),p_0)\bigr],\qquad
\det DB(\xi,s')=
-\frac{\partial}{\partial s'}
\bigl[\rho(q_0,\tau(\xi,s'))\bigr].
\]
The coordinate choices above are possible because
$|\partial_q\rho(z_0)|\ge\kappa$ and
$|\partial_p\rho(z_0)|\ge\kappa$:
each of these two-dimensional vectors has a component
of absolute value at least $\kappa/\sqrt2$.
The $C^2$ bound for $\rho$ gives a common neighborhood
of $(0,0)$ on which
$|\det DA|\ge\kappa/(2\sqrt2)$ and
$|\det DB|\ge\kappa/(2\sqrt2)$.
By the quantitative inverse function theorem
\cite[Proposition~2.5.6]{AbrahamMarsdenRatiu},
the preceding estimates give a common $r_0>0$
such that $A$ and $B$ admit $C^m$ inverse maps
on $V=(-r_0,r_0)^2$.
Define the manifold-valued coordinate parametrizations by
$\mathcal Q_{z_0}=\sigma\circ A^{-1}$ and
$\mathcal P_{z_0}=\tau\circ B^{-1}$.
These coordinate parametrizations satisfy
$\mathcal Q_{z_0}(0,0)=q_0$ and
$\mathcal P_{z_0}(0,0)=p_0$,
and their product image is contained in $U$.
The $C^m$ bound for $\rho$ gives uniform bounds
for the derivatives of $A$ and $B$ up to order $m$
on $V_Q$ and $V_P$, respectively.
The determinant lower bounds and the identities
$D(A^{-1})=(DA\circ A^{-1})^{-1}$ and
$D(B^{-1})=(DB\circ B^{-1})^{-1}$
therefore give uniform bounds for the first
derivatives of both inverse maps.
Repeated differentiation of
$A\circ A^{-1}=\operatorname{Id}$ and
$B\circ B^{-1}=\operatorname{Id}$ then gives,
by induction on the derivative order,
\begin{equation}\label{eq:normal-coordinate-bounds}
\begin{gathered}
\|D(A^{-1})\|_{C^{m-1}(V)}
+
\|D(B^{-1})\|_{C^{m-1}(V)}
\le C,\\
C^{-1}\le|\det D(A^{-1})(w)|\le C,
\qquad
C^{-1}\le|\det D(B^{-1})(w)|\le C
\quad(w\in V).
\end{gathered}
\end{equation}

For the remainder of the proof, we suppress the index $z_0$
and set $f(y,t,x,t')=
\rho\bigl(\mathcal Q(y,t),\mathcal P(x,t')\bigr)$.
The definitions of the coordinates give
$f(y,t,0,0)=t$ and $f(0,0,x,t')=-t'$.
Differentiating these identities gives
$f_y(0)=f_x(0)=f_{yy}(0)=f_{xx}(0)=0$,
$f_t(0)=1$, and $f_{t'}(0)=-1$.
Set $c=f_{yx}(0)$.
Substituting these derivatives into \eqref{eq:rot} gives
\[
\Rot(f)(0)
=
\det
\begin{pmatrix}
0&0&1\\
0&c&f_{xt}(0)\\
-1&f_{t'y}(0)&f_{t't}(0)
\end{pmatrix}
=c.
\]
Computing $\Rot(\rho)(z_0)$ in the coordinates
$(\eta,s,\xi,s')$, the chain rule gives
\[
c
=
\det D(A^{-1})(0)\,
\det D(B^{-1})(0)\,
\Rot(\rho)(z_0).
\]
Using \eqref{eq:normal-coordinate-bounds}
and the assumptions on $\rho$,
we obtain by the chain rule
\[
C^{-1}\le|c|\le C,
\qquad
\|f\|_{C^m(V\times V)}\le C,
\]
after enlarging $C$ if necessary.

To obtain the rescaled formula, write
$\zeta=(y,t,x,t')$ and define
\[
\delta_L\zeta=(Ly,L^2t,Lx,L^2t').
\]
Set $A_D=\max\{1,\sup_{\zeta\in D}|\zeta|_\infty\}$
and choose $0<L_0\le
\min\left\{1,\frac{r_0}{2A_D}\right\}$.
Then $\delta_L(D)\subset V\times V$
for $0\le L\le L_0$, so we can define
$G(L,\zeta)=f(\delta_L\zeta)$.

The derivative identities established above give
$G(0,\zeta)=\partial_LG(0,\zeta)=0$
and $\partial_L^2G(0,\zeta)=
2t f_t(0)+2t'f_{t'}(0)
+y^2f_{yy}(0)+2yx f_{yx}(0)+x^2f_{xx}(0)=2(t-t'+cyx)$.
Taylor's formula with integral remainder therefore gives
\[
G(L,\zeta)
=
L^2(t-t'+cyx)+L^3\Theta_L(\zeta),
\]
where
$\Theta_L(\zeta)=
\frac12\int_0^1(1-\theta)^2
(\partial_L^3G)(\theta L,\zeta)\,d\theta$.
Dividing by $L^2$ proves \eqref{eq:scaled-normal}.

It remains to bound the remainder uniformly in
$L$ and $z_0$. For $|\alpha|\le m-3$, the chain rule for
$\partial_\zeta^\alpha\partial_L^3G$
involves derivatives of $f$ of order at most
$|\alpha|+3\le m$.
Since $D$ is fixed and $L_0\le1$, it follows that
\[
\sup_{0\le L\le L_0}
\|(\partial_L^3G)(L,\cdot)\|_{C^{m-3}(D)}
\le C\|f\|_{C^m(V\times V)}
\le C.
\]
Differentiating the integral defining $\Theta_L$
therefore yields $\sup_{0<L\le L_0}
\|\Theta_L\|_{C^{m-3}(D)}
\le C$.
All constants are independent of $z_0$,
which completes the proof.
\end{proof}

To use \cref{lem:normal-form} in an operator estimate,
we need a partition of unity with bounds for the derivatives
after rescaling and for the number of sets containing each point.
The following lemma provides these bounds.
For the use of the scales $\delta$ and $\delta^2$
in local coordinates, see
\cite[Section~3, (3.21)--(3.22), pp.~139--140]
{PhongSteinRadonI}.
We write $\mathcal Q_\nu$ and $\mathcal P_\nu$
for the adapted coordinate parametrizations centered at
$z_\nu=(q_\nu,p_\nu)$, and set
$\delta_L(y,t,x,t')=(Ly,L^2t,Lx,L^2t')$.

The next lemma and \cref{thm:curved} use the same
curvature condition.

\noindent\textbf{Assumption.}
The functions $a\in C_c^\infty(Q\times P)$ and the real-valued
function $\rho\in C^2(Q\times P;\R)$ satisfy
\begin{equation}\label{eq:curved-nondegeneracy}
\Rot(\rho)(z)\ne0
\qquad\text{for every } z=(q,p)\in\supp a\cap\{\rho=0\}.
\end{equation}

\begin{lemma}\label{lem:plate-cover}


Let $r\ge0$ be an integer, let $c_*>0$, and fix $a\in C_c^\infty(Q\times P)$.
Let $\mathcal F$ be a compact subset of
$C^{r+3}(Q\times P;\mathbb R)$.\footnote{Throughout, we equip
$C^k(Q\times P;\mathbb R)$ with the topology of
uniform convergence of all coordinate derivatives
of order at most $k$ on every compact subset of
each coordinate neighborhood.}
Assume that every $\rho\in\mathcal F$ satisfies
\eqref{eq:curved-nondegeneracy} for this $a$.

There exist constants $A>0$, $C\ge1$, and $L_0>0$
such that, for every $\rho\in\mathcal F$
and $0<L\le L_0$, there are a finite index set $I$,
centers $z_\nu\in\{\rho=0\}$, and nonnegative functions
$b_\nu\in C_c^r(Q\times P)$ satisfying the following
properties.

\begin{enumerate}[label=\textup{(\roman*)}]

    \item
     $\sum_\nu b_\nu=1$ on
     $\supp a\cap\{|\rho|\le c_*L^2\}$ and
     $\supp b_\nu\subset Q_\nu\times P_\nu$,
     where $B_L=(-AL,AL)\times(-AL^2,AL^2)$,
     $Q_\nu=\mathcal Q_\nu(B_L)$, and
     $P_\nu=\mathcal P_\nu(B_L)$.

    \item
    The normal-form identity \eqref{eq:scaled-normal}
holds for $\mathcal Q_\nu,\mathcal P_\nu$
on $D=[-A,A]^4$, with
$C^{-1}\le|c_\nu|\le C$ and
$\|\Theta_{\nu,L}\|_{C^r(D)}\le C$.
The coordinate bounds \eqref{eq:normal-coordinate-bounds}
also hold with $m=r+3$.
All constants are independent of $\nu$ and $L$.
    Moreover,
    \[
    \bigl\|
    b_\nu\circ(\mathcal Q_\nu\times\mathcal P_\nu)
    \circ\delta_L
    \bigr\|_{C^r([-A,A]^4)}
    \le C.
    \]

    \item
    $\#\{\nu:q\in Q_\nu\}\le CL^{-1}$ and
    $\#\{\nu:p\in P_\nu\}\le CL^{-1}$
for every $q\in Q$ and $p\in P$.
\end{enumerate}
\end{lemma}

\begin{proof}

Set $K=\supp a$.
If $K=\varnothing$ or $\mathcal F=\varnothing$,
take $A=C=L_0=1$ and $I=\varnothing$.
Assume henceforth that both sets are nonempty.

\medskip
\noindent\textbf{Step 1. Local preparation near $\rho_0$.}\par\nobreak\smallskip

We first prove that the constants $A,C,L_0$ can be chosen uniformly
on a neighborhood of each fixed $\rho_0\in\mathcal F$.

Fix $\rho_0\in\mathcal F$ and set
$\Sigma_0=K\cap\{\rho_0=0\}$.
If $\Sigma_0=\varnothing$, put
$\delta=\min_K|\rho_0|>0$ and take
\[
\mathcal U_{\rho_0}
=\left\{\rho\in C^{r+3}(Q\times P;\mathbb R):
\|\rho-\rho_0\|_{C^0(K)}<\delta/2\right\}.
\]
Then $|\rho|>\delta/2$ on $K$ for
$\rho\in\mathcal U_{\rho_0}$.
Taking $c_*L_0^2\le\delta/4$ makes
$K\cap\{|\rho|\le c_*L^2\}$ empty for $0<L\le L_0$ and any $\rho\in\mathcal U_{\rho_0}$.
Thus $I=\varnothing$ satisfies the conclusions
uniformly on $\mathcal U_{\rho_0}$, with any fixed $A>0$.

Now assume that $\Sigma_0\ne\varnothing$.
By \eqref{eq:rot}, both $d_q\rho_0$ and $d_p\rho_0$
are nonzero on $\Sigma_0$.
By continuity and compactness, choose open neighborhoods
$\Sigma_0\subset U_0\Subset U\Subset W$
on which all three nondegeneracy conditions hold for $\rho_0$.
For each $z=(q_z,p_z)\in\Sigma_0$, choose local
product coordinates
$q=\sigma_z(\bar y,\bar t)$ and
$p=\tau_z(\bar x,\bar t')$
on an open neighborhood $O_z$ of $z$,
with $\sigma_z(0)=q_z$ and $\tau_z(0)=p_z$,
such that
$f_{0,z}:=\rho_0\circ(\sigma_z\times\tau_z)$
satisfies
$\partial_{\bar t}f_{0,z}(0)\ne0$ and
$\partial_{\bar t'}f_{0,z}(0)\ne0$.
By continuity, choose an open product coordinate box
$V_z\Subset U_0\cap O_z$ containing $z$
such that $\partial_{\bar t}f_{0,z}$ and
$\partial_{\bar t'}f_{0,z}$ are bounded away from zero
on $\widetilde V_z=(\sigma_z\times\tau_z)^{-1}(V_z)$.
The implicit function theorem gives a $C^{r+3}$ function
$\varphi_{0,z}$, defined on an open neighborhood of
$0\in\mathbb R^3$, such that
$\varphi_{0,z}(0)=0$ and
$f_{0,z}(w,\varphi_{0,z}(w))=0$.
After restricting its domain, we may assume that
its graph lies in $\widetilde V_z$.

Set $\pi(\bar y,\bar t,\bar x,\bar t')
=(\bar y,\bar t,\bar x)$.
Choose a compact product box $K_z\subset\mathbb R^3$
containing $\pi((\sigma_z\times\tau_z)^{-1}(z))=(0,0,0)$
in its interior, such that $\varphi_{0,z}$ is defined
on an open neighborhood of $K_z$ and
\[
\{(w,\varphi_{0,z}(w)):w\in K_z\}
\Subset\widetilde V_z.
\]
Choose an open product coordinate box
$U_z\Subset V_z$ containing $z$ and $\varepsilon_z>0$
such that, with
$\widetilde U_z=(\sigma_z\times\tau_z)^{-1}(U_z)$,
\[
\pi(\overline{\widetilde U_z})
+[-\varepsilon_z,\varepsilon_z]^3\subseteq K_z.
\]
Since $\Sigma_0$ is compact, choose $z_1,\ldots,z_N$
such that the boxes $U_{z_\ell}$ cover $\Sigma_0$.
Relabel these pairs as $(U_\ell,V_\ell)$.
Write $\sigma_\ell=\sigma_{z_\ell}$,
$\tau_\ell=\tau_{z_\ell}$, $K_\ell=K_{z_\ell}$, $O_\ell=O_{z_\ell}$
and put $f_{0,\ell}=f_{0,z_\ell}$ and
$\varphi_{0,\ell}=\varphi_{0,z_\ell}$.
Set $\varepsilon=\min_\ell\varepsilon_{z_\ell}>0$.
Write
$\widetilde U_\ell=(\sigma_\ell\times\tau_\ell)^{-1}(U_\ell)$
and
$\widetilde V_\ell=(\sigma_\ell\times\tau_\ell)^{-1}(V_\ell)$.
Compactness and the finiteness of the cover give
constants $B_0<\infty$ and $\kappa_1,\kappa>0$ such that
\[
\begin{gathered}
\|f_{0,\ell}\|_{C^{r+3}(\widetilde V_\ell)}\le B_0,\\
\min\{|\partial_{\bar t}f_{0,\ell}|,
|\partial_{\bar t'}f_{0,\ell}|\}\ge2\kappa_1,
\qquad
|\Rot(f_{0,\ell})|\ge2\kappa
\quad\text{on }\widetilde V_\ell.
\end{gathered}
\]

The compact graph inclusions allow us to choose $h>0$
such that, for every $\ell$,
\begin{equation}\label{eq:plate-graph-tube}
\{(w,s):w\in K_\ell,\,
|s-\varphi_{0,\ell}(w)|\le2h\}
\Subset\widetilde V_\ell.
\end{equation}
For $\rho\in C^{r+3}(Q\times P;\mathbb R)$,
write $f_{\rho,\ell}=\rho\circ(\sigma_\ell\times\tau_\ell)$.
For each $\ell$, define the continuous seminorm
\[
p_\ell(g)
=
\max_{|\alpha|\le r+3}
\sup_{u\in\overline{\widetilde V_\ell}}
\left|
\partial^\alpha
\bigl(g\circ(\sigma_\ell\times\tau_\ell)\bigr)(u)
\right|.
\]
Its continuity follows from the chosen topology,
since $\overline{V_\ell}$ is a compact subset of $O_\ell$.
For $\delta>0$, the set
\[
\mathcal U
=
\left\{
\rho\in C^{r+3}(Q\times P;\mathbb R):
p_\ell(\rho-\rho_0)<\delta
\text{ for every }\ell
\right\}
\]
is an open neighborhood of $\rho_0$, since there
are finitely many $\ell$.
For every $\rho\in\mathcal U$ and every $\ell$ we have 
$
\|f_{\rho,\ell}-f_{0,\ell}\|_{C^{r+3}(\widetilde V_\ell)}
\le p_\ell(\rho-\rho_0)<\delta.
$
Taking $\delta$ sufficiently small, the preceding bounds
for $f_{0,\ell}$ therefore give
\begin{equation}\label{eq:plate-local-bounds}
\begin{gathered}
\|f_{\rho,\ell}\|_{C^{r+3}(\widetilde V_\ell)}
\le2B_0,
\qquad
\|f_{\rho,\ell}-f_{0,\ell}\|_{C^0(\widetilde V_\ell)}
<\kappa_1h,\\
\min\{|\partial_{\bar t}f_{\rho,\ell}|,
|\partial_{\bar t'}f_{\rho,\ell}|\}\ge\kappa_1,
\qquad
|\Rot(f_{\rho,\ell})|\ge\kappa
\quad\text{on }\widetilde V_\ell.
\end{gathered}
\end{equation}

For $w\in K_\ell$, since
$f_{0,\ell}(w,\varphi_{0,\ell}(w))=0$ and
$|\partial_{\bar t'}f_{0,\ell}|\ge2\kappa_1$
on $\widetilde V_\ell$,
\eqref{eq:plate-graph-tube} and the mean value theorem
imply that the values
$f_{0,\ell}(w,\varphi_{0,\ell}(w)-h)$ and
$f_{0,\ell}(w,\varphi_{0,\ell}(w)+h)$
have opposite signs and absolute values at least $2\kappa_1h$.
By $\|f_{\rho,\ell}-f_{0,\ell}\|_{C^0(\widetilde V_\ell)}
<\kappa_1h$, the corresponding values for $f_{\rho,\ell}$
still have opposite signs.
We have already shown that $|\partial_{\bar t'}f_{\rho,\ell}|
\ge\kappa_1>0$ on $\widetilde V_\ell$.
By continuity and the connectedness of the product box
$\widetilde V_\ell$, the derivative
$\partial_{\bar t'}f_{\rho,\ell}$ has constant sign
on $\widetilde V_\ell$.
The intermediate value theorem
and strict monotonicity imply that the equation
$f_{\rho,\ell}(w,s)=0$ has a unique solution in
$(\varphi_{0,\ell}(w)-h,\varphi_{0,\ell}(w)+h)$.
Define $\varphi_{\rho,\ell}(w)$ to be this solution.
Then
\[
f_{\rho,\ell}(w,\varphi_{\rho,\ell}(w))=0,
\qquad
|\varphi_{\rho,\ell}(w)-\varphi_{0,\ell}(w)|<h.
\]
The implicit function theorem extends $\varphi_{\rho,\ell}$
to a $C^{r+3}$ function on an open neighborhood of $K_\ell$.
The bounds
$\|f_{\rho,\ell}\|_{C^{r+3}(\widetilde V_\ell)}\le2B_0$
and $|\partial_{\bar t'}f_{\rho,\ell}|\ge\kappa_1$
give uniform bounds for its derivatives of orders
$1,\ldots,r+3$ by repeated implicit differentiation.
Together with
$|\varphi_{\rho,\ell}-\varphi_{0,\ell}|<h$ on $K_\ell$,
this yields
\begin{equation}\label{eq:plate-graph-derivative-bound}
\|\varphi_{\rho,\ell}\|_{C^{r+3}(K_\ell)}
\le C_\varphi
\qquad(\rho\in\mathcal U,\ \ell=1,\ldots,N),
\end{equation}
where $C_\varphi$ is independent of $\rho$ and $\ell$.
For each $\ell$, we have
\begin{equation}\label{eq:plate-graph-inclusion}
\{(w,\varphi_{\rho,\ell}(w)):w\in K_\ell\}
\subseteq
\{(w,s):w\in K_\ell,\ |s-\varphi_{0,\ell}(w)|\le h\}
\Subset\widetilde V_\ell.
\end{equation}
Since there are finitely many $\ell$,
the distances of these graphs from
$\partial\widetilde V_\ell$ have a common positive
lower bound.

If $K\setminus\bigcup_\ell U_\ell\ne\varnothing$,
choose $\delta_0>0$ such that
$
3\delta_0<
\min_{z\in K\setminus\bigcup_\ell U_\ell}|\rho_0(z)|.
$
Such a choice is possible because this set is compact
and contains no zero of $\rho_0$.
If the set is empty, choose any $\delta_0>0$.
In either case,
$K\cap\{|\rho_0|\le3\delta_0\}\subseteq\bigcup_\ell U_\ell$.
Replace $\mathcal U$ by the open neighborhood
$\mathcal U\cap
\{\rho:\|\rho-\rho_0\|_{C^0(K)}<\delta_0\}$,
and set $\mathcal U_{\rho_0}=\mathcal U$.
Choose $0<L_1\le1$ with $c_*L_1^2\le\delta_0$.
Then
\[
K\cap\{|\rho|\le c_*L^2\}
\subseteq K\cap\{|\rho_0|\le2\delta_0\}
\subseteq\bigcup_\ell U_\ell
\qquad(\rho\in\mathcal U,\ 0<L\le L_1).
\]

\medskip
\noindent\textbf{Step 2. Grid points and adapted coordinates.}\par\nobreak\smallskip

Fix $\rho\in\mathcal U$ and $0<L\le L_1$, and set
$H=L^2$, $f_\ell=f_{\rho,\ell}$, and
$\varphi_\ell=\varphi_{\rho,\ell}$.
For fixed $\rho_0$, the charts, boxes, $\varepsilon$,
and the bounds obtained above are independent of
$\rho\in\mathcal U$ and $L$.
 We now construct a covering and
a partition of unity for $\rho$ with constants
independent of $\rho\in\mathcal U$ and $L$.

For the local construction in each box, we suppress
the index $\ell$ and write $\sigma,\tau,f,\varphi$
for $\sigma_\ell,\tau_\ell,f_\ell,\varphi_\ell$.
The identity
$f(\bar y,\bar t,\bar x,\varphi(\bar y,\bar t,\bar x))=0$
implies
$\partial_{\bar t}\varphi
=-(\partial_{\bar t}f)/(\partial_{\bar t'}f)$,
where the derivatives of $f$ are evaluated at
$\bar t'=\varphi(\bar y,\bar t,\bar x)$.
Since both $|\partial_{\bar t}f|$ and
$|\partial_{\bar t'}f|$ lie between $\kappa_1$
and $2B_0$, the preceding formula gives
\begin{equation}\label{eq:plate-graph-t-bound}
0<c\le|\partial_{\bar t}\varphi|\le C
\qquad\text{on }K_\ell,
\end{equation}
where $c,C$ depend only on
$B_0,\kappa_1$ and are independent of
$\rho\in\mathcal U$, $\ell$, and $L$.

Set $E_L=\supp a\cap\{|\rho|\le c_*H\}$.
Choose grids $y_i=iL$, $t_n=nH$, and $x_j=jL$
in the $\bar y,\bar t,\bar x$ variables, respectively,
for $i,j,n\in\mathbb Z$.
Define the finite index set
\[
I=\{(\ell,i,j,n):1\le\ell\le N,\ (y_i,t_n,x_j)\in K_\ell\}.
\]
For $\nu=(\ell,i,j,n)\in I$, write $\ell(\nu)=\ell$ and set
\[
s_\nu=\varphi_\ell(y_i,t_n,x_j),\qquad
q_\nu=\sigma_\ell(y_i,t_n),\qquad
p_\nu=\tau_\ell(x_j,s_\nu).
\]
For each $\ell$, the set
\[
\Gamma_\ell=
\bigl\{\bigl(\sigma_\ell(y,t),
\tau_\ell(x,\varphi_\ell(y,t,x))\bigr):(y,t,x)\in K_\ell\bigr\}
\]
is a compact subset of $V_\ell\cap\{\rho=0\}$.
Thus $\{z_\nu=(q_\nu,p_\nu):\nu\in I\}
\subseteq
\bigcup_{\ell=1}^N\Gamma_\ell
\subseteq U_0\cap\{\rho=0\}$.

Fix $\nu=(\ell,i,j,n)\in I$.
Write $\widetilde V_\ell=\Omega_\ell^Q\times\Omega_\ell^P$,
where $\Omega_\ell^Q,\Omega_\ell^P\subset\mathbb R^2$ are open product boxes.
Define $u_\nu(q)=\rho(q,p_\nu)$, $v_\nu(p)=-\rho(q_\nu,p)$,
and the coordinate maps
\begin{equation}\label{eq:plate-coordinate-maps}
A_\nu(\bar y,\bar t)=(\bar y-y_i,u_\nu(\sigma(\bar y,\bar t))),\qquad
B_\nu(\bar x,\bar t')=(\bar x-x_j,v_\nu(\tau(\bar x,\bar t'))).
\end{equation}
Here $A_\nu$ is defined on $\Omega_\ell^Q$ and
$B_\nu$ on $\Omega_\ell^P$.
Since $\rho(q_\nu,p_\nu)=0$, we have
$A_\nu(y_i,t_n)=B_\nu(x_j,s_\nu)=(0,0)$.
Differentiating \eqref{eq:plate-coordinate-maps}, we obtain
the Jacobian determinants
$\partial_{\bar t}f(\bar y,\bar t,x_j,s_\nu)$ and
$-\partial_{\bar t'}f(y_i,t_n,\bar x,\bar t')$, respectively.
By \eqref{eq:plate-local-bounds}, each determinant has
absolute value at least $\kappa_1$ and hence has constant
sign on its connected coordinate box.
By these derivative bounds, the second component of
$A_\nu$ is strictly monotone in $\bar t$ for fixed
$\bar y$, and that of $B_\nu$ is strictly monotone
in $\bar t'$ for fixed $\bar x$.
Since their first components are translations,
both maps are injective.
The uniform $C^2$ bounds for $A_\nu$ and $B_\nu$ from
\eqref{eq:plate-local-bounds}, their Jacobian lower bound
$\kappa_1$, and the common positive distance of the graphs
from the coordinate boundaries (see \eqref{eq:plate-graph-tube})
allow us to apply the quantitative inverse function argument
in the proof of \cref{lem:normal-form} uniformly for
$\rho\in\mathcal U$.
It gives a square $B_*=(-r_0,r_0)^2$, with $r_0>0$
independent of $\rho\in\mathcal U$, $\nu$, and $L$, such that
$B_*\subseteq A_\nu(\Omega_\ell^Q)\cap B_\nu(\Omega_\ell^P)$.
Applying the argument leading to \eqref{eq:normal-coordinate-bounds}
in the proof of \cref{lem:normal-form} to $A_\nu$ and $B_\nu$,
with the uniform $C^{r+3}$ bounds from \eqref{eq:plate-local-bounds}
and the Jacobian lower bound $\kappa_1$, gives the same estimates
for $A_\nu^{-1}$ and $B_\nu^{-1}$ on $B_*$ with $m=r+3$, uniformly in
$\rho\in\mathcal U$, $\nu$, and $L$.
Set $\mathcal Q_\nu=\sigma\circ A_\nu^{-1}$ and
$\mathcal P_\nu=\tau\circ B_\nu^{-1}$.
Write $(\eta,u)=A_\nu(\bar y,\bar t)$ and
$(\xi,v)=B_\nu(\bar x,\bar t')$; in particular,
$\eta=\bar y-y_i$ and $\xi=\bar x-x_j$.
We have $\mathcal Q_\nu(B_*)\times\mathcal P_\nu(B_*)\subset V_\ell$,
and for $(\eta,u),(\xi,v)\in B_*$,
\begin{equation}\label{eq:plate-coordinate-identities}
\rho(\mathcal Q_\nu(\eta,u),p_\nu)=u,
\qquad
\rho(q_\nu,\mathcal P_\nu(\xi,v))=-v.
\end{equation}

\medskip
\noindent\textbf{Step 3. Covering the sublevel set and obtaining the normal form.}\par\nobreak\smallskip

We now estimate the coordinate widths needed to cover $E_L$.
By \eqref{eq:plate-graph-derivative-bound}, choose $M\ge1$
such that $|D\varphi_\ell|\le M$ on every $K_\ell$.
Fix $R>M/(2c)$, where $c$ is the lower bound in
\eqref{eq:plate-graph-t-bound}. Both $M$ and $R$ can be
chosen independently of $\rho\in\mathcal U$.
For the covering argument, take $0<L\le1$
sufficiently small that
$E_L\subseteq\bigcup_\ell U_\ell$ and
$\max\{L/2,RL+H/2\}<\varepsilon$.
Fix $\ell\in\{1,\ldots,N\}$ and $(q,p)\in E_L\cap U_\ell$, with coordinates
$(\bar y,\bar t,\bar x,\bar t')$, and choose $y_i,x_j$ with
$|y_i-\bar y|,|x_j-\bar x|\le L/2$.
Since
$\pi(\overline{\widetilde U_\ell})+[-\varepsilon,\varepsilon]^3\subseteq K_\ell$,
we have
$(y_i,\lambda,x_j)\in K_\ell$ whenever $|\lambda-\bar t|\le RL$.
Moreover,
\[
|\varphi(y_i,\bar t,x_j)-\varphi(\bar y,\bar t,x_j)|\le ML/2.
\]
On the interval $[\bar t-RL,\bar t+RL]$, $\lambda\mapsto\varphi(y_i,\lambda,x_j)$
has derivative of constant sign and absolute value at least $c$.
The intermediate value theorem therefore gives a unique $t_*$ with
\begin{equation}\label{eq:plate-t-star}
\varphi(y_i,t_*,x_j)=\varphi(\bar y,\bar t,x_j),
\qquad |t_*-\bar t|\le ML/(2c)<RL.
\end{equation}
Choose $t_n$ with $|t_n-t_*|\le H/2$. The inequalities
$|t_n-\bar t|\le H/2+RL<\varepsilon$ and
$|y_i-\bar y|,|x_j-\bar x|\le L/2<\varepsilon$ imply that $(y_i,t_n,x_j)\in K_\ell$.
Thus $\nu=(\ell,i,j,n)\in I$. Using \eqref{eq:plate-t-star}
and the mean value theorem with the derivative bound from
\eqref{eq:plate-graph-derivative-bound}, we obtain
\begin{equation}\label{eq:plate-s-nu-bound}
|s_\nu-\varphi(\bar y,\bar t,x_j)|
\le M|t_n-t_*|\le MH/2.
\end{equation}
The point $\tau(x_j,\varphi(\bar y,\bar t,x_j))$ satisfies
$\rho(q,\tau(x_j,\varphi(\bar y,\bar t,x_j)))=0$.
The uniform upper bound for $|\partial_{\bar t'}f|$ on $\widetilde V_\ell$
from \eqref{eq:plate-local-bounds}, together with
\eqref{eq:plate-s-nu-bound}, gives $|u_\nu(q)|\le CH$.
Applying the mean value theorem in $\bar t'$ also gives
\begin{equation}\label{eq:plate-tprime-bound}
|\bar t'-\varphi(\bar y,\bar t,\bar x)|
\le\kappa_1^{-1}|\rho(q,p)|\le\kappa_1^{-1}c_*H.
\end{equation}
Since $H\le L$, the grid estimates give
$|\bar y-y_i|+|\bar t-t_n|\le(R+1)L$.
The triangle inequality, the bound in \eqref{eq:plate-tprime-bound},
the bound $|\varphi(\bar y,\bar t,\bar x)
-\varphi(\bar y,\bar t,x_j)|\le ML/2$,
and \eqref{eq:plate-s-nu-bound} give
$|\bar x-x_j|+|\bar t'-s_\nu|\le CL$.
Since $\widetilde V_\ell$ is a product box, integrating the mixed
second derivatives of $f$ gives
\begin{equation}\label{eq:plate-mixed-difference}
\begin{aligned}
&\bigl|\rho(q,p)-\rho(q,p_\nu)-\rho(q_\nu,p)
+\rho(q_\nu,p_\nu)\bigr|\\
&\qquad\le
C\bigl(|\bar y-y_i|+|\bar t-t_n|\bigr)
\bigl(|\bar x-x_j|+|\bar t'-s_\nu|\bigr)
\le CL^2=CH.
\end{aligned}
\end{equation}
Since $\rho(q_\nu,p_\nu)=0$, combining
\eqref{eq:plate-mixed-difference} with the bounds
$|\rho(q,p)|\le c_*H$ and $|u_\nu(q)|\le CH$
gives $|v_\nu(p)|\le CH$.
Thus there is a constant $a_0>0$, independent of $L$, such that
for every $\ell=1,\ldots,N$ and every
$(q,p)\in E_L\cap U_\ell$, there exist $i,j,n\in\mathbb Z$
with $\nu=(\ell,i,j,n)\in I$ and
\begin{equation}\label{eq:plate-grid-bounds}
|\bar y-y_i|,\ |\bar x-x_j|\le a_0L,
\qquad
|u_\nu(q)|,\ |v_\nu(p)|\le a_0H.
\end{equation}

Fix any $A>a_0$. The constants $C$ and $L_0$ in the statement
of the lemma can be chosen uniformly for $\rho\in\mathcal U$.
Choose $0<L_0\le\min\{1,L_1\}$ such that for $0<L\le L_0$,
$E_L\subset\bigcup_\ell U_\ell$ and
\[
\max\{L_0/2,RL_0+L_0^2/2\}<\varepsilon,
\qquad AL_0<r_0.
\]
With $B_L=(-AL,AL)\times(-AH,AH)$, we have $\overline{B_L}\subset B_*$.
Define $Q_\nu=\mathcal Q_\nu(B_L)$ and
$P_\nu=\mathcal P_\nu(B_L)$.
By \eqref{eq:plate-coordinate-maps} and the injectivity
of $A_\nu$ and $B_\nu$, we have
\begin{equation}\label{eq:plate-set-descriptions}
\begin{aligned}
Q_\nu&=\{q\in\sigma_\ell(\Omega_\ell^Q):
|\bar y-y_i|<AL,\ |u_\nu(q)|<AH\},\\
P_\nu&=\{p\in\tau_\ell(\Omega_\ell^P):
|\bar x-x_j|<AL,\ |v_\nu(p)|<AH\}.
\end{aligned}
\end{equation}
Since $A>a_0$, \eqref{eq:plate-grid-bounds} and
\eqref{eq:plate-set-descriptions} give
$E_L\subseteq\bigcup_{\nu\in I}(Q_\nu\times P_\nu)$.
For $D=[-A,A]^4$, the coordinate identities
\eqref{eq:plate-coordinate-identities} allow us to apply
the Taylor expansion from the proof of \cref{lem:normal-form}
to $\mathcal Q_\nu$ and $\mathcal P_\nu$, giving
\eqref{eq:scaled-normal}.
The bounds in \eqref{eq:plate-local-bounds} and
\eqref{eq:normal-coordinate-bounds} yield
$C^{-1}\le |c_\nu|\le C$ and
$\|\Theta_{\nu,L}\|_{C^r(D)}\le C$ uniformly in
$\rho\in\mathcal U$, $\nu$, and $L$.
This proves the normal-form and coordinate assertions in \textup{(ii)}.

\medskip
\noindent\textbf{Step 4. Overlap bounds.}\par\nobreak\smallskip

We verify \textup{(iii)} before constructing the partition.
Fix $\ell$, and write the compact product box as
$K_\ell=Y_\ell\times T_\ell\times X_\ell$,
where $Y_\ell,T_\ell,X_\ell$ are compact intervals.
We first count only the indices $\nu=(\ell,i,j,n)\in I$.
If $q\in Q_\nu$ for at least one such index, write
$q=\sigma(\bar y,\bar t)$ in this chart.
The inequality $|\bar y-y_i|<AL$ permits at most
$2A+2$ values of $i$, since the $y_i=iL$ have spacing $L$.
The condition $x_j=jL\in X_\ell$ permits at most $CL^{-1}$
values of $j$, since $X_\ell$ is a fixed bounded interval.
For each of these pairs with $y_i\in Y_\ell$ and $x_j\in X_\ell$,
define on $T_\ell$
\[
g_{i,j,q}(\lambda)
=f\bigl(\bar y,\bar t,x_j,\varphi(y_i,\lambda,x_j)\bigr).
\]
Since $q\in Q_\nu\subset\sigma_\ell(\Omega_\ell^Q)$,
we have $(\bar y,\bar t)\in\Omega_\ell^Q$.
For $\lambda\in T_\ell$, $(y_i,\lambda,x_j)\in K_\ell$,
so the graph inclusion \eqref{eq:plate-graph-inclusion} gives
$(x_j,\varphi(y_i,\lambda,x_j))\in\Omega_\ell^P$.
Thus $g_{i,j,q}$ is well defined on $T_\ell$.
Differentiating gives
\[
g_{i,j,q}'(\lambda)=
(\partial_{\bar t'}f)
\bigl(\bar y,\bar t,x_j,\varphi(y_i,\lambda,x_j)\bigr)
(\partial_{\bar t}\varphi)(y_i,\lambda,x_j).
\]
Both factors are continuous on the interval $T_\ell$, with
$|\partial_{\bar t'}f|\ge\kappa_1>0$ and
$|\partial_{\bar t}\varphi|\ge c>0$.
Thus $g_{i,j,q}'$ has constant sign on $T_\ell$ and
$|g_{i,j,q}'|\ge\kappa_1c$.
The set $\{\lambda\in T_\ell:|g_{i,j,q}(\lambda)|<AH\}$
is therefore an interval of length at most $2AH/(\kappa_1c)$.
As $u_\nu(q)=g_{i,j,q}(t_n)$ and the $t_n$ have spacing $H$,
at most $2A/(\kappa_1c)+2$ values of $n$ can occur.
Consequently, for this fixed $\ell$,
$\#\{\nu\in I:\ell(\nu)=\ell,\ q\in Q_\nu\}\le CL^{-1}$.
The analogous counting argument for fixed $p$, using
$|\partial_{\bar t}f|\ge\kappa_1$, gives
$\#\{\nu\in I:\ell(\nu)=\ell,\ p\in P_\nu\}\le CL^{-1}$.
Summing the bounds for fixed $q$ and fixed $p$ over
$\ell=1,\ldots,N$ proves \textup{(iii)}.

For the derivative estimates of the partition of unity,
we also need a uniform product overlap bound.
Fix $(q,p)\in Q\times P$ and $\ell$. We count the indices
$\nu=(\ell,i,j,n)\in I$ for which $(q,p)\in Q_\nu\times P_\nu$.
By \eqref{eq:plate-set-descriptions}, for each such index we have
$(q,p)\in V_\ell$ and the inequalities
$|\bar y-y_i|<AL$ and $|\bar x-x_j|<AL$ in this chart.
These inequalities permit at most $2A+2$ choices for each of $i$ and $j$.
The bound for $g_{i,j,q}$ then permits at most
$2A/(\kappa_1c)+2$ choices for $n$.
Summing over the same $N$ charts gives the product overlap bound
\begin{equation}\label{eq:plate-product-overlap}
\sup_{(q,p)\in Q\times P}
\#\{\nu\in I:(q,p)\in Q_\nu\times P_\nu\}\le C.
\end{equation}

\medskip
\noindent\textbf{Step 5. Partition of unity and derivative bounds.}\par\nobreak\smallskip

Choose a fixed nonnegative $\theta\in C_c^\infty((-A,A))$
with $\theta=1$ on $[-a_0,a_0]$. For $\nu=(\ell,i,j,n)\in I$, define
$\zeta_\nu$ on $\mathcal Q_\nu(B_*)\times\mathcal P_\nu(B_*)$ by
\begin{equation}\label{eq:cutoff-def}
\zeta_\nu(q,p)=
\theta\left(\frac{\bar y-y_i}{L}\right)
\theta\left(\frac{u_\nu(q)}H\right)
\theta\left(\frac{\bar x-x_j}{L}\right)
\theta\left(\frac{v_\nu(p)}H\right),
\end{equation}
and extend by zero to $Q\times P$.
Indeed, choose $b<A$ with $\supp\theta\subset[-b,b]$.
The support of $\zeta_\nu$ is contained in
\[
\mathcal Q_\nu([-bL,bL]\times[-bH,bH])
\times
\mathcal P_\nu([-bL,bL]\times[-bH,bH])
\Subset Q_\nu\times P_\nu.
\]
Since $\supp\theta\Subset(-A,A)$ and $\overline{B_L}\subset B_*$,
$\supp\zeta_\nu
\Subset\mathcal Q_\nu(B_*)\times\mathcal P_\nu(B_*)$.
The defining formula gives
$\zeta_\nu\in C^{r+3}
(\mathcal Q_\nu(B_*)\times\mathcal P_\nu(B_*))$.
Its zero extension therefore belongs to
$C_c^{r+3}(Q\times P)$.
For each $\ell$ and each $(q,p)\in E_L\cap U_\ell$,
the covering argument gives an index $\nu=(\ell,i,j,n)\in I$ such that
\[
\max\left\{
\frac{|\bar y-y_i|}{L},
\frac{|u_\nu(q)|}{H},
\frac{|\bar x-x_j|}{L},
\frac{|v_\nu(p)|}{H}
\right\}\le a_0.
\]
Since $\theta=1$ on $[-a_0,a_0]$, this index satisfies
$\zeta_\nu(q,p)=1$. For $\ell=1,\ldots,N$, set
\[
I_\ell=\{\mu\in I:\ell(\mu)=\ell\},
\qquad
S_\ell=\sum_{\mu\in I_\ell}\zeta_\mu.
\]
Then $S_\ell\ge1$ on $E_L\cap U_\ell$.

The compact set $K\cap\{|\rho_0|\le2\delta_0\}$
is contained in $\bigcup_\ell U_\ell$.
Choose fixed nonnegative functions $\chi_\ell\in C_c^\infty(U_\ell)$
such that $\sum_{\ell=1}^N\chi_\ell=1$ on this compact set.
These functions are independent of $\rho\in\mathcal U$ and $L$.
Since $E_{L_0}\subseteq K\cap\{|\rho_0|\le2\delta_0\}$
for every $\rho\in\mathcal U$, the same functions satisfy
$\sum_{\ell=1}^N\chi_\ell=1$ on every such $E_{L_0}$.
Choose also a fixed nonnegative $F\in C^\infty(\mathbb R)$ with
$F(s)=0$ for $s\le1/2$ and $F(s)=s^{-1}$ for $s\ge1$.
The function $F$ and all its derivatives are bounded on $\mathbb R$.
For $\nu\in I_\ell$, define
\[
b_\nu=\chi_\ell\zeta_\nu F(S_\ell).
\]
These functions are nonnegative, belong to $C_c^r(Q\times P)$,
and satisfy $\supp b_\nu\subset Q_\nu\times P_\nu$.
If $(q,p)\in E_L\subseteq E_{L_0}$ and $\chi_\ell(q,p)\ne0$,
then $(q,p)\in U_\ell$, so $S_\ell(q,p)\ge1$.
Consequently, on $E_L$,
\[
\sum_{\nu\in I}b_\nu
=\sum_{\ell=1}^N\chi_\ell S_\ell F(S_\ell)
=\sum_{\ell=1}^N\chi_\ell=1.
\]
This proves \textup{(i)}.

Fix $\nu=(\ell,i,j,n)\in I$.
We estimate
$\|S_\ell\circ(\mathcal Q_\nu\times\mathcal P_\nu)
\circ\delta_L\|_{C^r(D)}$.
Take $(q,p)\in Q_\nu\times P_\nu$ and write its adapted coordinates
as $(\eta,u)=(Ly,Ht)$ and $(\xi,v)=(Lx,Ht')$. Thus
\[
q=\mathcal Q_\nu(Ly,Ht),
\qquad
p=\mathcal P_\nu(Lx,Ht'),
\qquad
(y,t,x,t')\in\operatorname{int}D=(-A,A)^4.
\]
Suppose $\mu=(\ell,i',j',n')\in I_\ell$ and
$(q,p)\in\supp\zeta_\mu$.
For $\omega\in\{\mu,\nu\}$, the inclusion
$(q,p)\in Q_\omega\times P_\omega$ and
\eqref{eq:normal-coordinate-bounds} give
\[
|\sigma_\ell^{-1}(q)-\sigma_\ell^{-1}(q_\omega)|
+|\tau_\ell^{-1}(p)-\tau_\ell^{-1}(p_\omega)|\le CL.
\]
The triangle inequality therefore gives
\begin{equation}\label{eq:plate-center-distance}
|\sigma_\ell^{-1}(q_\mu)-\sigma_\ell^{-1}(q_\nu)|
+|\tau_\ell^{-1}(p_\mu)-\tau_\ell^{-1}(p_\nu)|\le CL.
\end{equation}
Define, for $(\eta,u),(\xi,v)\in B_L$,
\[
\begin{aligned}
\Delta_{\mu\nu}^Q(\eta,u)
&=\rho(\mathcal Q_\nu(\eta,u),p_\mu)
-\rho(\mathcal Q_\nu(\eta,u),p_\nu),\\
\Delta_{\mu\nu}^P(\xi,v)
&=\rho(q_\mu,\mathcal P_\nu(\xi,v))
-\rho(q_\nu,\mathcal P_\nu(\xi,v)).
\end{aligned}
\]
For $\omega\in\{\mu,\nu\}$, set
$\alpha_\omega=\tau_\ell^{-1}(p_\omega)$ and
$\beta_\omega=\sigma_\ell^{-1}(q_\omega)$.
For $0\le s\le1$, define
$\alpha_s=(1-s)\alpha_\nu+s\alpha_\mu$ and $\beta_s=(1-s)\beta_\nu+s\beta_\mu$.
Convexity gives $\alpha_s\in\Omega_\ell^P$ and
$\beta_s\in\Omega_\ell^Q$.
The fundamental theorem of calculus gives
\[
\begin{aligned}
\Delta_{\mu\nu}^Q(\eta,u)
&=\int_0^1
D_2f_\ell\bigl(A_\nu^{-1}(\eta,u),\alpha_s\bigr)
\cdot(\alpha_\mu-\alpha_\nu)\,ds,\\
\Delta_{\mu\nu}^P(\xi,v)
&=\int_0^1
D_1f_\ell\bigl(\beta_s,B_\nu^{-1}(\xi,v)\bigr)
\cdot(\beta_\mu-\beta_\nu)\,ds,
\end{aligned}
\]
where $D_1,D_2$ differentiate the first and second
$\mathbb R^2$ variables, respectively.
Differentiating under the integral sign and using the bounds in
\eqref{eq:plate-local-bounds}, \eqref{eq:normal-coordinate-bounds},
and \eqref{eq:plate-center-distance}, we obtain
\begin{equation}\label{eq:plate-phase-difference-bound}
\|\Delta_{\mu\nu}^Q\|_{C^r(B_L)}
+\|\Delta_{\mu\nu}^P\|_{C^r(B_L)}
\le C\bigl(
|\alpha_\mu-\alpha_\nu|
+|\beta_\mu-\beta_\nu|\bigr)
\le CL,
\end{equation}
where the derivatives are in $(\eta,u)$ and $(\xi,v)$, respectively.
By \eqref{eq:plate-coordinate-identities}, we have
\begin{equation}\label{eq:cutoff-phase-identities}
\frac{u_\mu(q)}H=t+\frac{\Delta_{\mu\nu}^Q(Ly,Ht)}H,
\qquad
\frac{v_\mu(p)}H=t'-\frac{\Delta_{\mu\nu}^P(Lx,Ht')}H.
\end{equation}
For nonnegative integers $a,b$ with $1\le a+b\le r$,
the chain rule and \eqref{eq:plate-phase-difference-bound} give
\[
\left|\partial_y^a\partial_t^b
\frac{\Delta_{\mu\nu}^Q(Ly,Ht)}H\right|
\le C\frac{L^aH^b}{H}L
=CL^{a+2b-1}\le C,
\]
and the same bound holds for
$\partial_x^a\partial_{t'}^b[\Delta_{\mu\nu}^P(Lx,Ht')/H]$.
Thus by \eqref{eq:cutoff-phase-identities},
\[
\left|
\partial_y^a\partial_t^b
\left(\frac{u_\mu(q)}{H}\right)
\right|
+
\left|
\partial_x^a\partial_{t'}^b
\left(\frac{v_\mu(p)}{H}\right)
\right|
\le C,
\qquad 1\le a+b\le r.
\]
Thus all positive-order derivatives through order $r$ of
$u_\mu(q)/H$ and $v_\mu(p)/H$ are uniformly bounded.
The arguments of the first and third factors in
\eqref{eq:cutoff-def}, with $\nu$ replaced by $\mu$, satisfy
\[
\frac{\bar y-y_{i'}}{L}=y+i-i',
\qquad
\frac{\bar x-x_{j'}}{L}=x+j-j'.
\]
Thus
$\zeta_\mu(\mathcal Q_\nu(Ly,Ht),\mathcal P_\nu(Lx,Ht'))
=\theta(y+i-i')
\theta(t+\Delta_{\mu\nu}^Q(Ly,Ht)/H)
\theta(x+j-j')
\theta(t'-\Delta_{\mu\nu}^P(Lx,Ht')/H)$.
Since $\theta$ and its derivatives are bounded,
the chain and product rules yield
$\left|\partial^\gamma[
\zeta_\mu\circ(\mathcal Q_\nu\times\mathcal P_\nu)
\circ\delta_L](y,t,x,t')\right|\le C$
for every $\gamma\in\mathbb N_0^4$ with $|\gamma|\le r$.
Set $\Phi_{\nu,L}=(\mathcal Q_\nu\times\mathcal P_\nu)\circ\delta_L$
and $z=(y,t,x,t')\in\operatorname{int}D$,
so that $\Phi_{\nu,L}(z)=(q,p)$.
If $(q,p)\notin\supp\zeta_\mu$, then
$\zeta_\mu\circ\Phi_{\nu,L}$ vanishes in a neighborhood of $z$.
By $\supp\zeta_\mu\subset Q_\mu\times P_\mu$
and the product overlap bound \eqref{eq:plate-product-overlap}, we have
$\#\{\mu\in I_\ell:(q,p)\in\supp\zeta_\mu\}\le C$.
Therefore, for every $\gamma\in\mathbb N_0^4$ with $|\gamma|\le r$,
\[
\left|\partial^\gamma(S_\ell\circ\Phi_{\nu,L})(z)\right|
\le
\sum_{\substack{\mu\in I_\ell\\
(q,p)\in\supp\zeta_\mu}}
\left|\partial^\gamma(\zeta_\mu\circ\Phi_{\nu,L})(z)\right|
\le C\,\#\{\mu\in I_\ell:(q,p)\in\supp\zeta_\mu\}.
\]
Applying the overlap bound and extending the estimate by continuity
to $\partial D$ gives $\bigl\|S_\ell\circ(\mathcal Q_\nu\times\mathcal P_\nu)
\circ\delta_L\bigr\|_{C^r(D)}\le C$.

In the same coordinates, $\zeta_\nu\circ(\mathcal Q_\nu\times\mathcal P_\nu)\circ\delta_L
=\theta(y)\theta(t)\theta(x)\theta(t')$.
The functions $\chi_\ell$ are fixed and smooth, so
\eqref{eq:normal-coordinate-bounds} and $L,H\le1$ give
$\bigl\|\chi_\ell\circ(\mathcal Q_\nu\times\mathcal P_\nu)
\circ\delta_L\bigr\|_{C^r(D)}\le C$.
Applying the product and chain rules to
$b_\nu=\chi_\ell\zeta_\nu F(S_\ell)$ now gives
\[
\bigl\|b_\nu\circ(\mathcal Q_\nu\times\mathcal P_\nu)
\circ\delta_L\bigr\|_{C^r(D)}
=\bigl\|(\chi_\ell\circ\Phi_{\nu,L})
(\zeta_\nu\circ\Phi_{\nu,L})
(F\circ S_\ell\circ\Phi_{\nu,L})\bigr\|_{C^r(D)}
\le C.
\]
This proves the remaining assertion of \textup{(ii)},
with a bound independent of $\rho\in\mathcal U$.

\medskip
\noindent\textbf{Step 6. Uniform constants on $\mathcal F$.}\par\nobreak\smallskip

For each $\rho_0\in\mathcal F$, we have constructed a neighborhood
$\mathcal U_{\rho_0}$ and a number $a_0(\rho_0)\ge0$
such that, for every fixed $A>a_0(\rho_0)$,
the construction gives \textup{(i)--(iii)} on
$\mathcal U_{\rho_0}$ with common constants $C,L_0$.
In the case $\Sigma_0=\varnothing$, take $a_0(\rho_0)=0$.
Both $\mathcal U_{\rho_0}$ and $a_0(\rho_0)$ are fixed
before $A$ is chosen.
By compactness, choose a finite subcover
$\mathcal F\subseteq\bigcup_{i=1}^s\mathcal U_{\rho_i}$.
Fix $A>\max_{1\le i\le s}a_0(\rho_i)$ and perform the construction
on each $\mathcal U_{\rho_i}$ with this same $A$.
Let $C_i\ge1$ and $L_{0,i}>0$ be the resulting constants.
Taking $C=\max_i C_i$ and $L_0=\min_i L_{0,i}$
proves the assertion for all $\rho\in\mathcal F$.

\end{proof}

The following lemma allows us to combine the local operator
estimates from the plate decomposition into a bound for
the full operator.

\begin{lemma}
\label{lem:support-schur}
Let $(P,\mu_P)$ and $(Q,\mu_Q)$ be measure spaces,
let $E$ be a Banach space, and let $I$ be finite.
For each $\nu\in I$, let $P_\nu\subset P$ and $Q_\nu\subset Q$
be measurable, and let
$T_\nu:L^2(P,\mu_P;E)\to L^2(Q,\mu_Q;E)$
be a bounded linear operator with
$T_\nu=1_{Q_\nu}T_\nu1_{P_\nu}$ and $\|T_\nu\|\le b$.
If $\sum_{\nu\in I}1_{P_\nu}\le N_P$ $\mu_P$-almost everywhere
and $\sum_{\nu\in I}1_{Q_\nu}\le N_Q$ $\mu_Q$-almost everywhere,
then
\[
\left\|\sum_{\nu\in I}T_\nu\right\|
\le b\sqrt{N_PN_Q}.
\]
\end{lemma}

\begin{proof}
Fix $f\in L^2(P,\mu_P;E)$.
Since $T_\nu=1_{Q_\nu}T_\nu$, we have
$(Tf)(q)=\sum_{\nu:q\in Q_\nu}(T_\nu f)(q)$
for almost every $q\in Q$.
The overlap bound on $Q$, the triangle inequality,
and Cauchy--Schwarz give $\|(Tf)(q)\|_E^2
\le N_Q\sum_{\nu\in I}\|(T_\nu f)(q)\|_E^2$.

Integrating over $Q$ and using
$T_\nu=T_\nu1_{P_\nu}$ and $\norm{T_\nu}\le b$,
we obtain
\[
\norm{Tf}_2^2
\le N_Q\sum_{\nu\in I}\norm{T_\nu f}_2^2
\le N_Qb^2\int_P
\left(\sum_{\nu\in I}1_{P_\nu}(p)\right)
\norm{f(p)}_E^2\,d\mu_P(p)
\le N_QN_Pb^2\norm f_2^2.
\]
Taking square roots and then the supremum
over $\norm f_2=1$ proves the claim.
\end{proof}

\subsection{Decay estimates for integral operators}

Fix $a\in C_c^\infty(Q\times P)$ and
$k\in\mathcal S(\mathbb R)$ with
$\supp\widehat k\Subset\mathbb R\setminus\{0\}$.
For a real-valued function $\rho$ and $R\ge1$, define
\[
(\mathcal K_Rf)(q)
=\int_P Rk(R\rho(q,p))a(q,p)f(p)\,d\mu_P(p).
\]
The following theorem extends the decay estimate
in \cref{thm:flat} to functions with nonzero
rotational curvature.

\begin{theorem}\label{thm:curved}
Let $\mathcal E$ be a class of type $>1$, and let
$\varepsilon>0$ be as in \cref{thm:flat}.
For every sufficiently large integer $m$, let
$\mathcal F\subset C^m(Q\times P;\R)$ be compact
and suppose that every $\rho\in\mathcal F$ satisfies
\eqref{eq:curved-nondegeneracy}.
Then there exists $C<\infty$ such that
\[
\|\mathcal K_R\otimes\Id_E\|_{L^2(P;E)\to L^2(Q;E)}
\le CR^{-\varepsilon/3}
\qquad(\rho\in\mathcal F,\ E\in\mathcal E,\ R\ge1).
\]

\end{theorem}

\begin{proof}

By the proof of \cref{cor:flat-amplitude}, we may take $r=8$.
Assume that $m\ge r+3$. The continuous inclusion
$C^m(Q\times P;\R)\hookrightarrow C^{r+3}(Q\times P;\R)$
shows that $\mathcal F$ is compact in $C^{r+3}(Q\times P;\R)$.
Since every $\rho\in\mathcal F$ satisfies
\eqref{eq:curved-nondegeneracy} for the fixed $a$,
we may apply \cref{lem:plate-cover} with this $r$ and $c_*=2$.
This gives constants $A,C,L_0$ with $L_0\le1$ that work
for all $\rho\in\mathcal F$.
Fix $\rho\in\mathcal F$.

For $1\le R\le R_0:=L_0^{-3}$, 
the bound
$|Rk(R\rho)a|\le R_0\|k\|_\infty|a|$
and the scalar Schur test
(see \cite[Theorem~5.2, p.~22]{HalmosSunder} and
\cite[Lemma~7.5, p.~52]{WolffHA}) give
$\|\mathcal K_R\otimes\Id_E\|
\le C_0R_0^{\varepsilon/3}R^{-\varepsilon/3}$,
where $C_0$ is independent of $\rho$, $R$, and $E$.
Assume henceforth that $R\ge R_0$, and set
\[
L=R^{-1/3},\qquad H=L^2,\qquad M=RH=R^{1/3}.
\]
Thus $L\le L_0$, $R=M^3$, and $ML=1$.
Choose $0\le\chi\le1$ in $C_c^\infty((-2,2))$ with
$\chi=1$ on $[-1,1]$, and write
$\mathcal K_R=\mathcal K_R^{\mathrm{near}}+
\mathcal K_R^{\mathrm{far}}$, where
$\mathcal K_R^{\mathrm{near}}$ has kernel
$Rk(R\rho)a\chi(\rho/H)$.
Since $\chi=1$ on $[-1,1]$, we have
$1-\chi(\rho(q,p)/H)=0$ whenever $|\rho(q,p)|\le H$.
Hence
\[
\supp\bigl(Rk(R\rho)a(1-\chi(\rho/H))\bigr)
\subseteq \supp a\cap\{|\rho|\ge H\}.
\]
On this support we have $R|\rho|\ge RH=M$.
Since $R=M^3$ and $0\le\chi\le1$, the Schwartz decay of $k$
of order $N+3$ gives
$\bigl|Rk(R\rho)a(1-\chi(\rho/H))\bigr|
\le C_N M^3(1+M)^{-(N+3)}
\le C_N M^{-N}$
for every integer $N\ge0$. The constant $C_N$ depends only on
$N$, $k$, and $\|a\|_\infty$, so the bound is uniform for
$\rho\in\mathcal F$.
Set
$K_Q=\operatorname{pr}_Q(\supp a)$ and
$K_P=\operatorname{pr}_P(\supp a)$.
Both sets are compact and have finite measure.
For $f\in L^2(P;E)$, the kernel bound gives
\[
\bigl\|[((\mathcal K_R^{\mathrm{far}})
\otimes\Id_E)f](q)\bigr\|_E
\le C_NM^{-N}1_{K_Q}(q)
\int_{K_P}\|f(p)\|_E\,d\mu_P(p)
\le C_NM^{-N}1_{K_Q}(q)
\mu_P(K_P)^{1/2}\|f\|.
\]
Taking the $L^2(Q)$ norm of the pointwise bound yields
\begin{equation}\label{eq:curved-far-bound}
\|(\mathcal K_R^{\mathrm{far}})\otimes\Id_E\|
_{L^2(P;E)\to L^2(Q;E)}\le C_N
\bigl[\mu_Q(K_Q)\mu_P(K_P)\bigr]^{1/2}M^{-N}
\le C_N'M^{-N}.
\end{equation}

For this $\rho$ and $L$, choose
$I$, $b_\nu$, $Q_\nu$, and $P_\nu$ by the construction
in the proof of \cref{lem:plate-cover}.
Retain the notation
$V_\ell,\sigma_\ell,\tau_\ell,
\Omega_\ell^Q,\Omega_\ell^P,A_\nu,B_\nu$
from that construction.
The finite subcover in the proof of \cref{lem:plate-cover}
gives a fixed finite collection of product charts
$\sigma_\ell\times\tau_\ell$ for the whole family $\mathcal F$.
On the corresponding compact coordinate boxes, smoothness
gives uniform $L^\infty$ bounds for the chart densities and
uniform $C^r$ bounds for $a\circ(\sigma_\ell\times\tau_\ell)$.
Since $\sum_\nu b_\nu=1$ on $\supp a\cap\{|\rho|\le2H\}$,
we have $\mathcal K_R^{\mathrm{near}}=\sum_{\nu\in I}T_\nu$, where
\[
(T_\nu f)(q)
=\int_P Rk(R\rho(q,p))a(q,p)b_\nu(q,p)
\chi(\rho(q,p)/H)f(p)\,d\mu_P(p).
\]
In particular, $T_\nu=1_{Q_\nu}T_\nu1_{P_\nu}$.
Put $D_0=(-A,A)^2$ and $D=[-A,A]^4$.
For $(y,t,x,t')\in D$, write
$q=\mathcal Q_\nu(Ly,Ht)$ and $p=\mathcal P_\nu(Lx,Ht')$.
By \cref{lem:plate-cover}\textup{(ii)},
\[
\begin{gathered}
\Psi_{\nu,L}(y,t,x,t'):=H^{-1}\rho(q,p)
=t-t'+c_\nu yx+L\Theta_{\nu,L}(y,t,x,t'),\\
C^{-1}\le|c_\nu|\le C,\quad
\|\Theta_{\nu,L}\|_{C^r(D)}\le C.
\end{gathered}
\]

In the chart $\ell=\ell(\nu)$, define the smooth positive
functions $w_\ell^Q,w_\ell^P$ by
$\sigma_\ell^*(d\mu_Q)
=w_\ell^Q(\bar y,\bar t)\,d\bar y\,d\bar t$
and
$\tau_\ell^*(d\mu_P)
=w_\ell^P(\bar x,\bar t')\,d\bar x\,d\bar t'$.
For $(y,t),(x,t')\in[-A,A]^2$, set
$j_\nu^Q(y,t)
=w_\ell^Q\bigl(A_\nu^{-1}(Ly,Ht)\bigr)
\left|\det D(A_\nu^{-1})(Ly,Ht)\right|$
and
$j_\nu^P(x,t')
=w_\ell^P\bigl(B_\nu^{-1}(Lx,Ht')\bigr)
\left|\det D(B_\nu^{-1})(Lx,Ht')\right|$.
Since $\mathcal Q_\nu=\sigma_\ell\circ A_\nu^{-1}$
and $\mathcal P_\nu=\tau_\ell\circ B_\nu^{-1}$,
the substitutions
$q=\mathcal Q_\nu(Ly,Ht)$ and
$p=\mathcal P_\nu(Lx,Ht')$ give
$d\mu_Q(q)=LHj_\nu^Q(y,t)\,dy\,dt$
and
$d\mu_P(p)=LHj_\nu^P(x,t')\,dx\,dt'$.
Since $w_\ell^Q,w_\ell^P$ are smooth on neighborhoods
of the compact sets $\overline{\Omega_\ell^Q}$ and
$\overline{\Omega_\ell^P}$, respectively, and there
are finitely many charts, we have
$
\max_{1\le\ell\le N}
\left\{
\|w_\ell^Q\|_{L^\infty(\Omega_\ell^Q)},
\|w_\ell^P\|_{L^\infty(\Omega_\ell^P)}
\right\}\le C_w<\infty,
$
where $C_w$ is independent of $\nu$ and $L$.
Together with the Jacobian bounds in \eqref{eq:normal-coordinate-bounds},
this gives $0<j_\nu^Q,j_\nu^P\le C$ on $\overline{D_0}$,
uniformly in $\rho\in\mathcal F$, $\nu$, and $L$.
Set
\[
\alpha_{\nu,L}(y,t,x,t')
=a(q,p)b_\nu(q,p)\chi(\Psi_{\nu,L}(y,t,x,t')).
\]
By \cref{lem:plate-cover}\textup{(ii)} and the product
and chain rules we have $\|\alpha_{\nu,L}\|_{C^r(D)}
\le C\|b_\nu\circ(\mathcal Q_\nu\times\mathcal P_\nu)
\circ\delta_L\|_{C^r(D)}
\le C$.
Moreover, $\supp b_\nu\Subset Q_\nu\times P_\nu$ gives
$\supp\alpha_{\nu,L}
\subseteq
\delta_L^{-1}\bigl(
(\mathcal Q_\nu\times\mathcal P_\nu)^{-1}(\supp b_\nu)
\bigr)
\Subset D_0\times D_0$.

For $f\in L^2(P;E)$, define $F(x,t')=j_\nu^P(x,t')f(\mathcal P_\nu(Lx,Ht'))$.
Since $0<j_\nu^P(x,t')\le C$ and
$d\mu_P(p)=LHj_\nu^P(x,t')\,dx\,dt'$
under $p=\mathcal P_\nu(Lx,Ht')$, we obtain
\[
\|F\|_{L^2(D_0;E)}^2
\le C\int_{D_0}j_\nu^P(x,t')
\|f(\mathcal P_\nu(Lx,Ht'))\|_E^2\,dx\,dt'
=\frac{C}{LH}\|f\|_{L^2(P_\nu;E)}^2
\le\frac{C}{LH}\|f\|_{L^2(P;E)}^2.
\]
Since $RLH=LM$, changing variables in $T_\nu$ gives
\begin{equation}\label{eq:scaled-block}
[(T_\nu\otimes\Id_E)f](\mathcal Q_\nu(Ly,Ht))=
L\iint_{D_0}M k(M\Psi_{\nu,L}(y,t,x,t'))
\alpha_{\nu,L}(y,t,x,t')F(x,t')\,dx\,dt'.
\end{equation}

Since $ML=1$, Fourier inversion gives
$M k(M\Psi_{\nu,L})
=\frac M{2\pi}\int_\R
e^{iM\omega(t-t'+c_\nu yx)}
e^{i\omega\Theta_{\nu,L}}\widehat k(\omega)\,d\omega$.
Define
\[
B_{\nu,L}(y,t,x,t',\omega)
=\alpha_{\nu,L}(y,t,x,t')
e^{i\omega\Theta_{\nu,L}(y,t,x,t')}\widehat k(\omega).
\]
Extending $B_{\nu,L}$ by zero for
$(y,t,x,t')\notin D$ gives
$B_{\nu,L}\in C_c^r(\mathbb R^5)$, with
$\supp B_{\nu,L}\subseteq D\times\supp\widehat k$ and 
$\|B_{\nu,L}\|_{C^r(\R^5)}\le C$.
Let $\mathsf A_{\nu,L}
=\mathsf A_{M,c_\nu,B_{\nu,L}}\otimes\Id_E$,
where $\mathsf A_{M,c,B}$ is defined in
\eqref{eq:flat-amplitude}.
By \cref{cor:flat-amplitude}, with
$D_x=D_y=\overline{D_0}=[-A,A]^2$
and a fixed compact annulus containing $\supp\widehat k$, we have 
$\|\mathsf A_{\nu,L}\|_{L^2(D_0;E)\to L^2(D_0;E)}
\le CM^{-\varepsilon}$.
The right-hand side of \eqref{eq:scaled-block} is
$L(\mathsf A_{\nu,L}F)(y,t)$.
Since $d\mu_Q(q)=LHj_\nu^Q(y,t)\,dy\,dt$
under $q=\mathcal Q_\nu(Ly,Ht)$,
$0<j_\nu^Q(y,t)\le C$, and
$\|F\|_{L^2(D_0;E)}^2
\le \frac{C}{LH}\|f\|_{L^2(P;E)}^2$,
we obtain
\[
\|(T_\nu\otimes\Id_E)f\|_{L^2(Q;E)}^2\le C LH\,L^2
\|\mathsf A_{\nu,L}F\|_{L^2(D_0;E)}^2\le C L^2M^{-2\varepsilon}\|f\|_{L^2(P;E)}^2.
\]
Thus $\|T_\nu\otimes\Id_E\|\le CLM^{-\varepsilon}$,
uniformly in $\rho\in\mathcal F$, $\nu$, $L$, and $E$.
By \cref{lem:plate-cover}\textup{(iii)},
$N_Q:=\sup_{q\in Q}\#\{\nu\in I:q\in Q_\nu\}\le CL^{-1}$
and
$N_P:=\sup_{p\in P}\#\{\nu\in I:p\in P_\nu\}\le CL^{-1}$.
If $I=\varnothing$, then $\mathcal K_R^{\mathrm{near}}=0$.
Otherwise, \cref{lem:support-schur} and the bound
for each $T_\nu\otimes\Id_E$ give
\[
\|\mathcal K_R^{\mathrm{near}}\otimes\Id_E\|
\le CLM^{-\varepsilon}\sqrt{N_PN_Q}
\le C'M^{-\varepsilon}.
\]
Choosing an integer $N>\varepsilon$ in
\eqref{eq:curved-far-bound} and combining the near and far
estimates gives $\|\mathcal K_R\otimes\Id_E\|
\le CM^{-\varepsilon}=CR^{-\varepsilon/3}$ for $R\ge R_0$.
Together with the bound for $1\le R\le R_0$, this proves
the estimate uniformly for $\rho\in\mathcal F$ and
$E\in\mathcal E$.
\end{proof}

\section{A Hölder estimate for spherical averages}
\label{sec:sphere}

Throughout this section, let $d\sigma$ denote surface measure
on $S^2$. For the circle-averaging operators introduced
in Section~2, we use the notation
\[
(T_\delta f)(x)
=(T_\delta^{S^2}f)(x)
=\int_{\{y\in S^2:\ip{x}{y}=\delta\}}
f(y)\,d\sigma_{x,\delta}(y),
\qquad x\in S^2,
\]
for $\delta\in(-1,1)$ and $f\in C(S^2)$.
At $\delta=1$ and $\delta=-1$, we set
$T_1f(x)=f(x)$ and $T_{-1}f(x)=f(-x)$, respectively.
To prove the estimate for $T_\delta-T_0$
required in \cref{thm:criterion}, we first verify the rotational-curvature condition
in \cref{thm:curved} for the functions
$\rho_\delta(x,y)=\ip{x}{y}-\delta$.

\begin{lemma}\label{lem:sphere-curvature}
For $\delta\in(-1,1)$, let
$\rho_\delta(x,y)=\ip{x}{y}-\delta$ on $S^2\times S^2$.
Then $\Rot(\rho_\delta)(x,y)\ne0$ whenever $\ip{x}{y}=\delta$.
\end{lemma}

\begin{proof}
Fix $(x,y)\in S^2\times S^2$ with $\ip{x}{y}=\delta$,
and set $s=\sqrt{1-\delta^2}>0$.
Define
\[
a=\frac{y-\delta x}{s},
\qquad
b=\frac{\delta y-x}{s},
\qquad
n=\frac{x\times y}{s}.
\]
Then $(a,n)$ and $(b,n)$ are bases of
$T_xS^2$ and $T_yS^2$, respectively.
For $u\in T_xS^2$ and $v\in T_yS^2$,
\[
d_x\rho_\delta[u]=\ip{u}{y},
\qquad
d_y\rho_\delta[v]=\ip{x}{v},
\qquad
d_y d_x\rho_\delta[u,v]=\ip{u}{v}.
\]
Consequently, in these bases the matrix in
\eqref{eq:rot} is
\[
\begin{pmatrix}
0&s&0\\
-s&\delta&0\\
0&0&1
\end{pmatrix},
\]
whose determinant is $s^2=1-\delta^2>0$.
The coordinate invariance of nonvanishing proves the assertion.
\end{proof}

We will approximate the circle average $T_\delta$ using
smooth kernels in the variable $\ip{x}{y}-\delta$.
The choice below gives an initial smoothing term and
successive differences whose Fourier transforms are supported
in annuli away from zero, as required by \cref{thm:curved}.
Choose a real-valued even $\kappa\in\cS(\R)$ such that
$\widehat\kappa=1$ on $[-1,1]$ and
$\supp\widehat\kappa\subset[-2,2]$,
and set $k(s)=2\kappa(2s)-\kappa(s)$.
Then $\int_{\R}\kappa(s)\,ds=1$,
$\widehat k(\omega)=\widehat\kappa(\omega/2)-\widehat\kappa(\omega)$,
and $\supp\widehat k\subset\{\omega:1\le|\omega|\le4\}$.
Define the initial smoothing operator $\mathsf S_\delta$ and,
for $j\ge0$, the difference operators $\mathsf D_{\delta,j}$ by
\begin{align*}
 \mathsf S_\delta f(x)
 &=\frac1{2\pi}\int_{S^2}
 \kappa(\ip{x}{y}-\delta)f(y)\,d\sigma(y),\\
 \mathsf D_{\delta,j}f(x)
 &=\frac1{2\pi}\int_{S^2}
 2^jk\bigl(2^j(\ip{x}{y}-\delta)\bigr)f(y)\,d\sigma(y),
 \quad j\ge0.
\end{align*}

\begin{proposition}\label{prop:spherical-annular}
Let $\mathcal E$ be a class of type $>1$.  There exist $C<\infty$ and
$\eta>0$ such that for all $E\in\mathcal E$, $j\ge0$, and
$\abs{\delta}\le1/2$,
\[
 \norm{\mathsf D_{\delta,j}\otimes\Id_E}_{L^2(S^2;E)\to L^2(S^2;E)}
 \le C2^{-j\eta}.
\]
\end{proposition}

\begin{proof}

Fix $m$ sufficiently large for \cref{thm:curved} and set
$\mathcal F=\{\rho_\delta:|\delta|\le1/2\}$.
Since $\rho_\delta-\rho_{\delta'}=\delta'-\delta$ is constant,
the map $\delta\mapsto\rho_\delta$ is continuous from
$[-1/2,1/2]$ into $C^m(S^2\times S^2;\R)$.
Thus $\mathcal F$ is compact.

Set $a\equiv(2\pi)^{-1}$ on $S^2\times S^2$.
By \cref{lem:sphere-curvature}, every $\rho\in\mathcal F$
satisfies \eqref{eq:curved-nondegeneracy} for this $a$.
Apply \cref{thm:curved} to $\mathcal F$ with $Q=P=S^2$
and the choices of $a$ and $k$ above.
For $\rho=\rho_\delta$ and $R=2^j$, we have
$\mathcal K_R=\mathsf D_{\delta,j}$.
The desired estimate follows with $\eta=\varepsilon/3$.

\end{proof}

To estimate $T_\delta-T_0$, we also need bounds for
$\mathsf D_{\delta,j}-\mathsf D_{0,j}$ and
$\mathsf S_\delta-\mathsf S_0$.

\begin{lemma}\label{lem:regular-difference}
There exists a constant $C>0$ such that, for every
Banach space $E$, every integer $j\ge0$, and every
$\delta\in[-1/2,1/2]$,
\[
\begin{aligned}
\|(\mathsf D_{\delta,j}-\mathsf D_{0,j})\otimes\Id_E\|
&\le C\min\{1,2^j|\delta|\},\\
\|(\mathsf S_\delta-\mathsf S_0)\otimes\Id_E\|
&\le C|\delta|.
\end{aligned}
\]
\end{lemma}

\begin{proof}
For $h\in L^1([-1,1])$ and $f\in C(S^2)$,
spherical coordinates give
\begin{equation}\label{eq:spherical-disintegration}
\frac1{2\pi}\int_{S^2}h(\ip{x}{y})f(y)\,d\sigma(y)
=\int_{-1}^1 h(s)T_sf(x)\,ds.
\end{equation}
Applying \eqref{eq:spherical-disintegration} with $f=1$,
we obtain, for every $x\in S^2$,
\[
\begin{aligned}
&\frac1{2\pi}\int_{S^2}
\left|2^jk\bigl(2^j(\ip{x}{y}-\delta)\bigr)
      -2^jk\bigl(2^j\ip{x}{y}\bigr)\right|\,d\sigma(y)\\
&\quad=\int_{-1}^1
\left|2^jk\bigl(2^j(s-\delta)\bigr)-2^jk(2^js)\right|\,ds\\
&\quad\le\min\{2\|k\|_{L^1(\R)},\,
                  2^j|\delta|\|k'\|_{L^1(\R)}\}
\le C\min\{1,2^j|\delta|\}.
\end{aligned}
\]
By symmetry in $x,y$, integration with respect to $d\sigma(x)$
gives the same bound for every fixed $y\in S^2$.
Let $K$ denote the kernel of
$\mathsf D_{\delta,j}-\mathsf D_{0,j}$.
For $f\in L^2(S^2;E)$, the Bochner integral inequality gives
\[
\left\|\int_{S^2}K(x,y)f(y)\,d\sigma(y)\right\|_E
\le
\int_{S^2}|K(x,y)|\,\|f(y)\|_E\,d\sigma(y).
\]
Applying the scalar Schur test to $|K|$, we obtain
$\|(\mathsf D_{\delta,j}-\mathsf D_{0,j})\otimes\Id_E\|
\le C\min\{1,2^j|\delta|\}$.
Replacing $2^jk(2^j\cdot)$ by $\kappa$ in this argument yields
$\|(\mathsf S_\delta-\mathsf S_0)\otimes\Id_E\|\le C|\delta|$.
\end{proof}

With these preparations, we can now prove the desired
estimate for $T_\delta-T_0$.

\begin{theorem}
\label{thm:spherical-holder}
Let $\mathcal E$ be a class of type $>1$.
There exist constants $C>0$ and $\alpha\in(0,1/2]$ such that
for every $E\in\mathcal E$ and $\delta\in[-1,1]$,
\[
 \norm{(T_\delta-T_0)\otimes\Id_E}_{L^2(S^2;E)\to L^2(S^2;E)}
 \le C\abs{\delta}^\alpha.
\]
\end{theorem}

\begin{proof}
The case $\delta=0$ is immediate.
The circle-averaging operators satisfy
$\|T_\delta\otimes\Id_E\|\le1$ on $L^2(S^2;E)$
for every Banach space $E$ and every $\delta\in[-1,1]$;
see \cite[Lemme~2.2 and the preceding paragraph]{LafforgueStrongT}
and \cite[Section~2.10, the discussion immediately following (3)]
{deLaSalleSL3}.
Hence, for $1/2<|\delta|\le1$ and $0<\alpha\le1/2$,
\[
\|(T_\delta-T_0)\otimes\Id_E\|
\le2\le2\sqrt2\,|\delta|^\alpha.
\]
Fix $0<|\delta|\le1/2$ and $E\in\mathcal E$.
For $f\in C(S^2)\otimes E$, parametrizing the circles by
$y=sx+\sqrt{1-s^2}\,e$, $e\in S^2\cap x^\perp$,
and using uniform continuity of $f$ gives
\[
\lim_{s\to\tau}
\|((T_s-T_\tau)\otimes\Id_E)f\|_{C(S^2;E)}=0
\qquad(\tau\in\{0,\delta\}).
\]
For an integer $N\ge0$, set
$\kappa_N(s)=2^{N+1}\kappa(2^{N+1}s)$.
Telescoping gives the finite identity
\begin{equation}\label{eq:spherical-telescoping}
\kappa(s)+\sum_{j=0}^N2^jk(2^js)=\kappa_N(s)
\qquad(s\in\R).
\end{equation}
For $\tau\in\{0,\delta\}$,
\eqref{eq:spherical-telescoping} and
\eqref{eq:spherical-disintegration} give
\[
\left[\left(\mathsf S_\tau+\sum_{j=0}^N\mathsf D_{\tau,j}\right)
\otimes\Id_E\right]f
=\int_{-1}^1\kappa_N(s-\tau)(T_s\otimes\Id_E)f\,ds.
\]
Since $\kappa\in\cS(\R)$ and $\int_{\R}\kappa=1$,
the functions $\kappa_N$ form an approximate identity.
Moreover, $|\tau|\le1/2$ implies
\[
\left|\int_{-1}^1\kappa_N(s-\tau)\,ds-1\right|
\le\int_{|u|\ge1/2}|\kappa_N(u)|\,du
\longrightarrow0.
\]
The continuity at $\tau$ established above, together with
$\sup_{s\in[-1,1]}\|(T_s\otimes\Id_E)f\|_{C(S^2;E)}
\le\|f\|_{C(S^2;E)}$, therefore gives
\[
\begin{aligned}
&\left\|
\left[\left(\mathsf S_\tau+\sum_{j=0}^N\mathsf D_{\tau,j}\right)
\otimes\Id_E\right]f
-(T_\tau\otimes\Id_E)f
\right\|_C
=\left\|
\int_{-1}^1\kappa_N(s-\tau)(T_s\otimes\Id_E)f\,ds
-(T_\tau\otimes\Id_E)f
\right\|_C\\
&\le
\int_{-1}^1|\kappa_N(s-\tau)|
\|((T_s-T_\tau)\otimes\Id_E)f\|_C\,ds
+
\left|\int_{-1}^1\kappa_N(s-\tau)\,ds-1\right|
\|(T_\tau\otimes\Id_E)f\|_C\\
&\le
\|\kappa\|_{L^1(\R)}
\sup_{|s-\tau|<h}
\|((T_s-T_\tau)\otimes\Id_E)f\|_C
+3\|f\|_{C(S^2;E)}\int_{|u|\ge h}|\kappa_N(u)|\,du.
\end{aligned}
\]
For any $\varepsilon>0$, choose $h$ small and then $N$ large
so that both terms are less than $\varepsilon/2$.
This  proves convergence in $C(S^2;E)$.
For integers $N\ge J\ge0$,
\cref{lem:regular-difference,prop:spherical-annular} give
\[
\begin{aligned}
&\left\|\left[\left(
\mathsf S_\delta-\mathsf S_0
+\sum_{j=0}^N(\mathsf D_{\delta,j}-\mathsf D_{0,j})
\right)\otimes\Id_E\right]f\right\|_{L^2(S^2;E)}\\
&\quad\le C\left(
|\delta|+|\delta|\sum_{j=0}^J2^j
+\sum_{j=J+1}^N2^{-j\eta}
\right)\|f\|_{L^2(S^2;E)}\\
&\quad\le C\bigl(|\delta|2^J+2^{-J\eta}\bigr)
\|f\|_{L^2(S^2;E)}.
\end{aligned}
\]
Since $\|g\|_{L^2(S^2;E)}
\le\sigma(S^2)^{1/2}\|g\|_{C(S^2;E)}$,
subtracting the limits for $\tau=\delta$ and $\tau=0$ gives
\[
\left[\left(
\mathsf S_\delta-\mathsf S_0
+\sum_{j=0}^N(\mathsf D_{\delta,j}-\mathsf D_{0,j})
\right)\otimes\Id_E\right]f\longrightarrow
((T_\delta-T_0)\otimes\Id_E)f
\quad\text{in }L^2(S^2;E).
\]
Letting $N\to\infty$ in the finite-sum estimate therefore gives
\[
\|((T_\delta-T_0)\otimes\Id_E)f\|_{L^2(S^2;E)}
\le C\bigl(|\delta|2^J+2^{-J\eta}\bigr)
\|f\|_{L^2(S^2;E)}.
\]
By density, this estimate holds for all $f\in L^2(S^2;E)$.
Choose an integer $J\ge0$ such that
$2^J\le|\delta|^{-1/(1+\eta)}<2^{J+1}$.
Then $|\delta|2^J+2^{-J\eta}
\le(1+2^\eta)|\delta|^{\eta/(1+\eta)}$
and hence
\[
 \norm{(T_\delta-T_0)\otimes\Id_E}
 \le C\abs{\delta}^{\eta/(1+\eta)}.
\]
Choosing $\alpha=\min\{1/2,\eta/(1+\eta)\}$ completes the proof.
\end{proof}

\section{\texorpdfstring{Strong Banach property $(T)$}{Strong Banach property (T)}}
\label{sec:strongT}

\begin{proof}[Proof of \cref{thm:main}]
Since $\mathcal E$ is of type $>1$, \cref{thm:spherical-holder}
gives constants $C<\infty$ and $\alpha\in(0,1/2]$ such that
\[
\|(T_\delta-T_0)\otimes\Id_E\|_{L^2(S^2;E)\to L^2(S^2;E)}
\le C|\delta|^\alpha
\qquad(E\in\mathcal E,\ \delta\in[-1,1]).
\]
Applying \cref{thm:criterion} proves the theorem.
\end{proof}

\begin{remark}\label{rem:quantitative}
In \cref{thm:main}, the exponent $t$ in
\cref{def:strong-banach-T} can be chosen uniformly over
spaces with common type bounds.
More precisely, for every $p\in(1,2]$, every $T<\infty$, 
and every length function $\ell$ on $\mathrm{SL}_3(\mathbb R)$, 
there exists $t=t(p,T,\ell)>0$ such that the same exponent $t$ may be
used for every space $E$ satisfying $T_p(E)\le T$.
Indeed, applying \cref{thm:spherical-holder} to
$\{E:T_p(E)\le T\}$ gives constants $C$ and $\alpha$
depending only on $p,T$.
The assertion then follows from \cref{rem:criterion-uniform}.
\end{remark}

\begin{proof}[Proof of \cref{cor:lattice-strongT}]
Let $G=\mathrm{SL}_3(\R)$ and let $\Gamma$ be a lattice in $G$.
Since $\mathcal E$ is of type $>1$, \cref{thm:spherical-holder}
gives constants $C>0$ and $\alpha\in(0,1/2]$ such that
\[
\|(T_\delta-T_0)\otimes\Id_E\|_{L^2(S^2;E)\to L^2(S^2;E)}
\le C|\delta|^\alpha
\qquad(E\in\mathcal E,\ \delta\in[-1,1]).
\]
We apply the second case of \cite[Theorem~1.5]{deLaSalleLattices} with $n=2$.
Indeed, $G$ is connected and simple, with Lie algebra
$\mathfrak{g}=\mathfrak{sl}_3(\mathbb{R})$, and
$T_\delta^{(2)}=T_\delta$ in our notation.
The preceding estimate therefore verifies condition~(1.5)
of \cite{deLaSalleLattices} with $\beta=2\alpha\in(0,1]$.
Hence $\Gamma$ has property~$(\ast_{\mathcal E})$ as defined in
\cite[Section~1, (1.1)--(1.2)]{deLaSalleLattices}.

Since $\mathcal E$ was an arbitrary class of type $>1$,
it follows from \cite[Corollary~3.5]{deLaSalleKazhdan}
and the equivalence established in
\cite[Section~2.9]{deLaSalleSL3}
that $\Gamma$ has strong Banach property~$(T)$
in the sense of \cref{def:strong-banach-T},
with the growth exponent uniform on $\mathcal E$.

\end{proof}

\section*{Acknowledgments}
We acknowledge the use of OpenAI's ChatGPT in developing
the arguments and take full responsibility for the
mathematical content.

\end{document}